\documentclass[%
	a4paper,
	twoside,
	10pt,
    toc=bibliography
]{article}%

    \usepackage[ngerman, english]{babel} 
    \usepackage{amsmath} 
    \usepackage{extarrows} 
    \usepackage{amsthm}    
    \usepackage[arrow, matrix, curve]{xy}   
    \usepackage{amssymb} 
    \usepackage[usenames, dvipsnames]{xcolor}
    \usepackage{minibox} 
    \usepackage{leftidx} 
    \usepackage{MnSymbol} 
    \usepackage{fancyhdr}			
	\usepackage{faktor}
	\usepackage{graphicx} 
	\usepackage{mathrsfs} 
	\usepackage[refpage, noprefix]{nomencl}
	\usepackage{makeidx}
	
	\usepackage[style=alphabetic, backend=bibtex]{biblatex}
    \bibliography{ms.bib}

\title{The Carlitz Logarithm as a Period Morphism for Local $G$-Shtukas}
\author{Paul Breutmann}
\date{}

	\newcommand\klg{\leqslant}   
	\newcommand\grg{\geqslant}    
	\newcommand\skl{\preccurlyeq}    
    \newcommand\ds{\displaystyle}
	\newcommand\N{{\mathbb N}} 
	\newcommand\Z{{\mathbb Z}} 
	\newcommand\R{{\mathbb R}}
	\newcommand\C{{\mathbb C}}
	\newcommand\M{{\mathbb M}}

	\newcommand\F{{\mathbb F}}

	\newcommand\PP{{\mathbb P}}
	\newcommand\q{{\mathfrak q}}
	\newcommand\p{{\mathfrak p}}
	\newcommand\oo{\mathcal{O}}
	\newcommand\gs{{RZ_{b,\mu}}}
	\newcommand\gsn{{RZ_{b,\mu,n}}}
    \newcommand\gf{\mathscr{RZ}_{b,\mu}}
    \newcommand\Qf{\widetilde{\mathcal{Q}}}

    \newcommand\Qs{\underline{\mathcal{Q}}}
    \newcommand\Qsr{{\mathcal{Q}}^{rig}}        
    \newcommand\flfts{{\mathbf{fSch}_{flft}}}

	\newcommand\bo{\hfill$\Box$}
    
    \newcommand\sortkey[1]{}
    \newcommand\sem[1]{\lsem #1 \rsem}
    \newcommand\semm{{\lsem z,z^{-1}\}}}
    \newcommand\semmm{{\lsem z,z^{-1}\}\left[\frac{1}{\lper}\right]}}

    \newcommand\Hom{\mbox{Hom}}
    \newcommand\fzet{\F_q\sem{\zeta}}
    \newcommand\fzett{\F_q\semr{\zeta}}
    \newcommand\nilp{{Nilp_{\F_q\sem{\zeta}}}}
    \newcommand\lper{l_{-}}
    \newcommand{\semr}[1]{(\mkern-4.5mu(#1)\mkern-4.5mu)}    
    
    \newcommand{\grk}[1]{\langle#1 \rangle}
    \newcommand\m{M}
    \newcommand\n{N}
    
    \newcommand\uklr{{{U^{rig}_{k,l}}}}
    \newcommand\op[1]{{#1}^{op}}
    \newcommand{\widesim}{{\scalebox{2}[1]{$\sim$}}}
    \newcommand\vecc[2]{(\begin{smallmatrix}
#1\\#2
\end{smallmatrix})}
    \newcommand\bvecc[2]{\begin{pmatrix}
#1\\#2
\end{pmatrix}}
\newcommand\casecomp{$\underline\M=\big(\F_q\lsem z\rsem ^2,
(\begin{smallmatrix}
z&0\\
0&1
\end{smallmatrix}
)\big)$\ }

	 \newlength\breite
	 \newcommand{\bruch}[3][0.5pt]{
	 \setbox0=\hbox{\ensuremath{#2}}
	 \setbox1=\hbox{\ensuremath{\diagup}}
 	 \setbox2=\hbox{\ensuremath{#3}}
	 \setlength\breite{1.2\ht0}
	 \addtolength{\breite}{1.2\ht2}
  	 \raisebox{0.5\ht0}{\usebox0}
 	 \mkern-10mu
  	\rotatebox{45}
	{\rule[-0.8\ht2]{\breite}{0.5pt}}
  	\mkern-10mu
  	\raisebox{-0.8\ht2}{\usebox2}\ 
}
    \newcommand\isom{\xrightarrow{\,\smash{\raisebox{-0.60ex}    
               {\ensuremath{\scriptstyle\sim}}}\ }}
    \newcommand\llmapsto{\longleftarrow\hspace{-4.3pt}
    \shortmid\ }
	 
	 \newlength\proofeinzug
	\newcommand\prof[1]{
   		\textsl{Proof:} #1 \bo
	} 
	 
	\newcommand\df[2]{
	        \begin{definition}
            [#1]
            #2 
            \end{definition}
            }
	\newcommand\rem[2]{\begin{remark}
            [#1]#2
            \end{remark}
            }
    \newcommand\lem[2]{\begin{lemma}
            [#1]#2
            \end{lemma}
            }
    \newcommand\ko[2]{\begin{korollar}
            [#1]#2
            \end{korollar}
            }
    \newcommand\prop[2]{\begin{pr}
            [#1]#2
            \end{pr}
            }
    \newcommand\theo[2]{\begin{theorem}
            [#1]#2
            \end{theorem}
            }
     \newcommand\bsp[2]{\begin{bssp}
            [#1]#2
            \end{bssp}
            }
            
    \newcommand\ger[1]{}

\usepackage[left=1.7cm, right=1.9cm,top=1cm,bottom=1.6cm, includeheadfoot]{geometry}
	
		\makeatletter
		\renewcommand{\section}{%
    		\@startsection {section}{1}{\z@}%
         		{-3.5ex plus -1ex minus -.2ex}%
          	{2.3ex plus.2ex}%
                 	{\normalfont\Large\bfseries\centerline}}
	\makeatother
		\makeatletter
\renewcommand\subsubsection{\@startsection{subsubsection}{3}{\z@}{-3.25ex
    plus -1ex minus -.2ex}{1.5ex plus .2ex}{\bfseries\centerline}}
	\makeatother

	\newtheorem{satz}{Satz}[section]
	 \newtheorem{definition}[satz]{Definition}
	 \newtheorem{lemma}[satz]{Lemma}
	 \newtheorem{korollar}[satz]{Corollary}
	 \newtheorem{pr}[satz]{Proposition}
	 \newtheorem{theorem}[satz]{Theorem}
	 \theoremstyle{definition}
	 \newtheorem{bssp}[satz]{Example}
	 \newtheorem{remark}{Remark}

\begin{document}






\pagenumbering{arabic}

\maketitle

\begin{abstract}
\noindent Local shtukas are the function field analogs for $p$-divisible groups. Similar to the $p$-adic theory, one defines Rapoport-Zink functors and Rapoport-Zink spaces for these local shtukas. 
The associated Hodge-Pink structures are described uniquely by a morphism,
called the period morphism of the moduli problem. 
We will prove the ind-representability of the Rapoport-Zink functor in a particular case and compute the corresponding Rapoport-Zink space as well as the corresponding period morphism. In this case, the period morphism is given by the Carlitz logarithm.
\end{abstract}
\mbox{}\\
\mbox{}
\tableofcontents
\mbox{}\\
\mbox{}

\section{Introduction}

Let $G$ be a split, reductive group over $\F_q$ and $T\subset G$ a split maximal Torus. Hartl and Viehmann define in \cite{HV11} local $G$-shtukas for these groups as well as the Rapoport-Zink functor that classifies local $G$-shtukas bounded by a dominant cocaracter $\mu\in X_{\star}(T)$ together with a quasi-isogeny to a fixed local $G$-shtuka $(\M,b)$ over $\F_q$. They proved that this functor is ind-representable by a closed ind-subscheme of the affine Grassmanian $\widehat{Gr}$. This Rapoport-Zink space is locally formally of finite type over $\F_q\sem{\zeta}$ and its reduced subscheme is the affine Deligne-Lusztig variety $X_{\mu}(b)$, for which there is a dimension formula. In \cite{HV16} Hartl and Viehmann define period spaces of weakly admissible Hodge-Pink structures on $\M$. One can associate with every point in the Rapoport-Zink space such a Hodge-Pink structure. In \cite[6]{HV16} Hartl and Viehmann prove that this is  uniquely described by a period morphism from the generic fiber of the Rapoport-Zink space to the period space. 
In this article we will compute an explicit example for this theory, where $G=Gl_2$, $(\M, \tau_{\M}) =\Big(\F_q\lsem z\rsem ^2,\big(
\begin{smallmatrix}
z&0\\
0&1
\end{smallmatrix}\big)
\Big)$
and $\mu=(1,0)$. In this case, the period morphism is given by the Carlitz logarithm. 

We fix a local field of characteristic $p$, which has residue field $\F_q$. After choosing a uniformizer $z$, we can identify this field with $\F_q\semr{z}$ and its ring of integers with $\F_q\sem{z}$. We define $\nilp$ as the category of schemes $S$ over $Spec\ \F_q\sem{\zeta}$ such that $\zeta$ is locally nilpotent in $\oo_S$. 
Let $S$ be a scheme in $\nilp$.
We denote by $\oo_S\sem{z}$ the sheaf of $\oo_S$-algebras, which associates with every open affine subscheme $U=Spec\ A\subset S$ the $A$-algebra $\oo_S(U)\sem{z}=A\sem{z}=A\hat\otimes_{\F_q}\F_q\sem{z}$.  We define $\oo_S\sem{z}[\frac{1}{z-\zeta}]$ to be the sheaf associated with the presheaf $U\mapsto \oo_S(U)\sem{z}[\frac{1}{z-\zeta}]$. Note that since $\zeta$ is nilpotent, we have $\oo_S\sem{z}[\frac{1}{z-\zeta}]=\oo_S\sem{z}[\frac{1}{z}]$.
\\
Let $M$ be a $\oo_S\sem{z}$-modules on $S$ and let $g:S'\to S$ be a morphism of schemes in $\nilp$, then $g^{*} M:=g^{-1}M\otimes_{g^{-1}(\oo_S\sem{z})}\oo_S'\sem{z}$ is naturally a $\oo_{S'}\sem{z}$-module. If no confusion is expected, we denote it also by $M\otimes_{\oo_S\sem{z}}\oo_{S'}\sem{z}$ or even by $M_{S'}$.
There is a special case of this pullback, namely where $g$ is equal to the absolute $q$ Frobenius $\sigma:S\to S$. We denote by $\sigma$ this Frobenius as well as the induced endomorphism $\sigma:\oo_S\sem{z}\to \oo_S\sem{z}$, which is defined by $\sigma(z)=z$ and $\sigma(a)=a^q$ for sections $a$ in $\oo_S$. Therefore, we have
\nomenclature{\sortkey{oM}$\sigma^{*}\m$}{$\sigma^{*}\m:= \m\otimes_{\oo_S\sem{z},\sigma}\oo_S\sem{z}$ for some $\oo_S\sem{z}$-module $M$}
$\sigma^{*}\m:= \m\otimes_{\oo_S\sem{z},\sigma}\oo_S\sem{z}$.
Given a morphism $f:\m\to \n$ of $\oo_S\sem{z}$-modules we often write 
\nomenclature{\sortkey{of}$\sigma^{*} f$}{denotes the morphism $f\otimes id:\sigma^{*}\m\to\sigma^{*}\n$ for some $f:M\to N$}
$\sigma^{*} f$ for the induced morphism $f\otimes id:\sigma^{*}\m\to\sigma^{*}N$. 
In addition we define $\m[\frac{1}{z-\zeta}]:= \m\otimes_{\oo_S\sem{z}}\oo_S\sem{z}[\frac{1}{z-\zeta}]$ and we also write $f$ for the induced morphism $f\otimes id:\m[\frac{1}{z-\zeta}]\to \n[\frac{1}{z-\zeta}]$. Note again that $\m[\frac{1}{z-\zeta}]$ equals $\m[\frac{1}{z}]:= \m\otimes_{\oo_S\sem{z}}\oo_S\sem{z}[\frac{1}{z}]$ as long as $\zeta$ is nilpotent.\\

Now a local shtuka of rank $r$ over $S$ is a pair $\underline{\m}=(\m, \tau_{\m})$, which consists of a locally free $\oo_S\sem{z}$-module $\m$ of rank $r$ and an isomorphism $\tau_{\m}:\sigma^{*}\m[\frac{1}{z-\zeta}]\isom \m[\frac{1}{z-\zeta}]$.
A morphism $f:(\m,\tau_{\m})\to (\n,\tau_{\n})$ of local shtukas over S is a morphism of $\oo_S\sem{z}$-modules $f:\m\to \n$ which satisfies $\tau_{\n}\circ \sigma^{*} f=f\circ\tau_{\m}$.
A quasi-isogeny $f:(\m,\tau_{\m})\to (\n,\tau_{\n})$ of local shtukas over $S$ is an isomorphism $f:\m[\frac{1}{z}]\to \n[\frac{1}{z}]$ of $\oo_S\sem{z}\sem{\frac{1}{z}}$-modules, which satisfies $\tau_{\n}\circ \sigma^{*} f=f\circ\tau_{\m}$.

For every morphism $g:S'\to S$ in $\nilp$ the pullback of $\underline M$ defined by $\underline M_{S'}:=(g^{*} M,\tau_M\otimes id)$ 
is naturally a local shtuka over $S'$. 
We remark that the category of the defined local shtukas of rank $r$ is equivalent to the category of local $GL_r$-shtukas as in \cite[\S 4]{HV11}.

We have the following notion of a bounded local shtuka, which is used to define the Rapoport-Zink functor. It will also be crucial for the Hodge-Pink structures associated with local shtukas. 
Let $\mu_1\grg\dots\grg\mu_r$ be a decreasing sequence of integers. A local shtuka $(\m,\tau_{\m})$ of rank r over $S\in \nilp$ is bounded by $\mu=(\mu_1,\dots,\mu_r)$ if\linebreak 
\begin{equation}\label{neubound}
\bigwedge\nolimits^{\!i}\tau_{\m}(\bigwedge\nolimits^{\!i}\sigma^{*}\m)\subseteq(z-\zeta)^{\mu_{r-i+1}+\dots+\mu_r}\cdot\bigwedge\nolimits^{\!i}\m\qquad \mbox{\ for }1\klg i\klg r,
\end{equation}
where the inclusion has to be an equality for $i=r$ (see \cite{HV11} Def.\ 4.3).
With this notation we can now define the Rapoport-Zink functor and formulate the maintheorem of this article.

We fix a local shtuka $\underline\M=(\M,\tau_\M)$ of rank $r$ over $\F_q$ and a $\mu=(\mu_1,\dots, \mu_r)\in \Z^r$ with $\mu_1\grg \dots\grg\mu_r$. After choosing a basis for the module $\M$, $\tau_\M$ is given by a matrix $b\in GL_r\big(\F_q\semr{z}\big)$.

\df{}{\label{rfunc}\mbox{}\\ The Rapoport-Zink functor $\gf$ is defined by
\begin{eqnarray*}
    \gf:\quad\op{(\nilp)}    \quad\longrightarrow & Set\\
    S\quad\longmapsto&\left\{\
    (\m,\tau_{\m}, \eta)\ {\Bigg\vert}\ 
        \minibox{$(\m,\tau_M)$ is a local shtuka
        over $S$ \\ bounded by $\mu$ and $\eta$ is a quasi-isogeny \\$\eta: M[\frac{1}{z{-\zeta}}]\isom \M_S[\frac{1}{z{-\zeta}}]$}
        \right\}/\sim\\
    (f:S'\to S)\quad \longmapsto &
        \gf(f):\gf(S)\to \gf(S') \\
        &(\m,\tau_{\m}, \eta)\mapsto (\m_{S'},\tau_{\m}\otimes id, \eta\otimes id)
\end{eqnarray*} 
where by definition $(\m,\tau_M,\eta_N)$ is isomorphic to $(\n,\tau_N,\eta_N)$ if and only if there exists an isomorphism\linebreak $g:(M,\tau_M)\to (\n,\tau_N)$ with $\eta_M=\eta_N\circ g$. We also write $(\m,\tau_M,\eta_N)\sim (\n,\tau_N,\eta_N)$. Note that $g$ is unique with this condition. Also note that by rigidity of quasi-isogenies \cite[Prop 3.9]{HV11} the quasi-isogeny $\eta$ is uniquely determined by its restriction to $V(\zeta)\subset S$.
}
In the second section we will prove the main theorem of this article, which is the following:

\theo{}{\label{monsterproposition}
Let $(\M, \tau_{\M}) :=(\M, b) :=\Big(\F_q\lsem z\rsem ^2,\big(
\begin{smallmatrix}
z&0\\
0&1
\end{smallmatrix}\big)
\Big)$
and $\mu=(1,0)$. Then the Rapoport-Zink functor $\gf$ is ind-representable by the formal scheme $\ds\gs:=\coprod_{(i,j)\in \Z^2}\ Spf\ \F_q\lsem \zeta,h\rsem$. 
If we write $\ds\gs:=\coprod_{(i,j)\in \Z^2}\ U_{ij}$ with $U_{ij}=\ds\varinjlim_n U_{ij}^n=\varinjlim_n Spec\ \F_q[\zeta,h]/(\zeta,h)^{n}$ then the universal ind object $\ds x^{univ}=(x_n^{univ})_{n\in\N}\in \varprojlim_n \gf(\coprod_{i,j\in\Z^2} U_{ij}^n)$ is given by
$$
x_n^{univ}\Big|_{U_{ij}^n}=\left(\F_q\sem{\zeta,h}/(\zeta, h)^{q^n}\sem{z}^2,\ \begin{pmatrix}
z-\zeta & 0 \\ h & 1 
\end{pmatrix}
,\ 
\begin{pmatrix}
z^i\prod_{i=0}^{n-1}\frac{z}{z-\zeta^{q^i}} & 0 \\
-z^j\sum_{i=0}^{n-1}\frac{{h}^{q^i}}{(z-\zeta)\dots(z-\zeta^{q^i})}& z^j
\end{pmatrix}\right)
$$
}

The goal of the third section is to associate with every point $x\in \gf(S)$ a Hodge Pink structure over $S$ on the associated z-isocrystal of $\underline\M$. There is a representable functor $\Qf$ that classifies these Hodge-Pink structures. We will therefore define maps $\gamma(S):\gf(S)\to \Qf(S^{rig})$. In the fourth section we will see that these maps define our period morphism, which is given by the Carlitz logarithm.\\

\textbf{Acknowledgements.} I would like to thank my advisor Urs Hartl, for all his helpful discussions. During the work of this project, the author was supported by the SFB 878 "Groups, Geometry \& Actions" of the German Science Foundation (DFG), the CNRS and the ERC Advanced Grant 742608 "GeoLocLang". 
\section{A Particular Rapoport-Zink Space}\label{sectionrapoportzinkspace}

\label{particular}
\subsubsection{Proof of Theorem \ref{monsterproposition}}

Let $\gsn=\coprod_{\Z^2}\ Spec\ \F_q[\zeta,h]/(\zeta,h)^n$ so that $\varinjlim_n \gsn=\gs$.
The ind-representability of the functor $\gf$ by $\gsn$ means by definition that this functor is isomorphic to the functor $\varinjlim_n Hom(-,\gsn)$. 
We know by \cite[Theorem 6.3]{HV11} that the Rapoport-Zink functor $\gf$ is representable by some formal scheme locally formally of finite type over $\F_q\sem{z}$. 
Therefore the functor $\gf$ is already defined by its restriction to affine noetherian connected schemes \mbox{$Spec\ R\in \nilp$}.
We will describe explicitly the isomorphism of functors $nat:Hom(-,\gs)\to \gf$. It suffices to define $nat(S)$ for all affine noetherian connected schemes \mbox{$Spec\ R\in \nilp$} and we will write $nat(R)$ instead of $nat(Spec\ R)$. We denote by $I:=\sqrt{(0)}$ the nilradical of $R$ which implies that 
$R/I$ is a reduced $\F_q\sem{z}$-algebra that has no non-trivial idempotents.\\
We write $R\semr{z}:=R\sem{z}[\frac{1}{z}]$, which equals $R\lsem z\rsem [\frac{1}{z-\zeta}]=R\semr{z}$ as long as $\zeta$ is nilpotent in $R$.
After choosing Zariski locally a basis for $\m$ the quasi-isogeny $\eta$ and the morphism $\tau_M$ can be viewed as matrices in $Gl_2\left(R\semr{z}\right)$ satisfying	
$\tau_M=\eta^{-1}b\sigma^{*} \eta$. The boundedness of the local  shtuka $\underline M$ by $\mu=(1,0)$ means by definition that $\det\tau_M\in(z-\zeta)R\lsem z\rsem ^{*}$ and that $\underline M$ is effective (i.e. $\tau_M(\sigma^{*} M)\subset M$) and, consequently, $\tau_M\in Mat_2(R\lsem z\rsem)$. 
Recall that the affine Grassmanian $Gr$ for $Gl_2$ is the sheafification of the presheaf $A\mapsto Gl_2(A\semr{z})/Gl_2(A\sem{z})$ for $\F_q$ algebras $A$. It is an ind-scheme over $\F_q$. We set $\widehat{Gr}:=Gr\times_{\F_q}Spf\ \F_q\sem{\zeta}$ and define:

\begin{equation}\label{deflambdamenge}
\Lambda(R):=\left\{\overline g\in \widehat{Gr}(R) \Big|\quad\minibox{locally we have $ g^{-1}b\sigma^{*} g\in 
Mat_2(R\sem{z})
$\\ $det(g^{-1}b\sigma^{*} g) \in (z-\zeta)R\lsem z\rsem ^{*}$}\right\}
\end{equation}
where $g$ denotes a representative of the equivalence class $\overline g$. 
The conditions on $\overline g$ do not depend on the representative and $\Lambda$ is functorial in $R$.
We claim that we have a bijection
\begin{align}
    \gf(Spec\ R) &\xleftrightarrow{1:1} \Lambda(R) \label{bilam} \\
    (M, \tau_M, \eta)&\longmapsto\ \ \overline{\eta}\nonumber \\
    \left(R\lsem z\rsem ^2, g^{-1}b\sigma^{*} g, g\right)&\llmapsto  \overline g \nonumber\quad.
\end{align}

The boundedness condition on $(M, \tau_M, \eta)$ implies $\overline\eta\in\Lambda(R)$. For $(M, \tau_M, \eta_N)\sim (N, \tau_N, \eta_N)$ an isomorphism is given by $h\in Gl_2(R\sem{z})$ with $\eta_M=\eta_N\circ h$ and $h\circ \tau_M=\sigma^\star \tau_N\circ h$. Hence both maps in \eqref{bilam} are well-defined and mutually inverse. We define ${J}$ as the group of quasi-isogenies $f:\underline\M_{\F_q^{alg}}\to \underline\M_{\F_q^{alg}}$ also denoted by $QIsog(\underline\M_{\F_q^{alg}})$. We have:
\begin{align*}
&{J}=\left\{h\in Gl_2(\F_q^{alg}\semr{z}):hb=b\sigma^{*} h\right\}\\
&\quad=\left\{\big(\begin{smallmatrix}h_{11}&h_{12}\\ h_{21}&h_{22}\end{smallmatrix}\big)\in Gl_2\big(\F_q^{alg}\semr{z}\big):\big(\begin{smallmatrix}zh_{11}&h_{12}\\ zh_{21}&h_{22}\end{smallmatrix}\big)=\big(\begin{smallmatrix}z\sigma(h_{11})&z\sigma(h_{12})\\\sigma(h_{21})&\sigma(h_{22})\end{smallmatrix}\big)\right\}=\begin{pmatrix}\F_q\semr{z}^{*} &0\\ 0&\F_q\semr{z}^{*}\end{pmatrix}\end{align*}
We remark that ${J}$ are the $\F_q$-valued points of an algebraic group $J_b$ which can be defined for all $b$ as in \cite{Vie06} and which operates on the connected components of the affine Deligne-Lusztig variety.
We define an operation of ${J}$ on $\Lambda(R)$ by
\begin{equation}\label{joperation}
{J}\times \Lambda(R)\to \Lambda(R),\quad (h,\overline g)\mapsto \overline{hg}.
\end{equation}
It will help us to determine the sets $\Lambda(R)$ and hence $\gf(Spec\ R)$. 
The element $\overline{hg}$ lies indeed in $\Lambda(R)$ since we have $(hg)^{-1}b\sigma^{*}(hg)=g^{-1}h^{-1}b\sigma^{*} h\sigma^{*} g=g^{-1}b\sigma^{*} g$ and $\overline{g}\in \Lambda(R)$. The operation is equivariant under the maps $\Lambda(R)\to \Lambda(R')$ from above. By the bijection \eqref{bilam} this operation clearly transfers to an operation on $\gf(Spec\ R)$.
We will use the following lemma \ref{transitivoperation} to determine at first the $R/I$-valued points of $\gf$. These points will already give us interesting information about the formal scheme $\gs$.\\

\lem{}{\label{transitivoperation} 
Let $Spec\ R$ be a connected noetherian scheme in $\nilp$ and $I=\sqrt{(0)}$.
The group
${J}
$ operates transitively on $\Lambda(R/I)$, i.e. for all $\overline g\in \Lambda(R/I)$ there exists an element $h\in {J}
$ with $(h,\overline g)\mapsto \overline{id}$ or equivalently $\overline{h^{-1}}=\overline g$. 
In particular we have $\Lambda(R/I)\simeq {J}/Stab(id)$.}

\prof{
We first prove the assertion in the case that $R/I$ is an integral domain. In this case we have $ord_z(a^{-1})=-ord_z(a)$ for $a\in R/I\semr{z}^{{*}}$ and $ord_z(a\cdot a')=ord_z(a)+ord_z(a')$ for $a,a'\in R/I\semr{z}$.
Let $g=\Big(\begin{smallmatrix}
\tilde\alpha&\tilde\beta\\
\tilde\gamma&\tilde\delta
\end{smallmatrix}\Big)\in \Lambda(R/I)$ with $\tilde\alpha, \tilde\beta, \tilde\gamma, \tilde\delta\in R/I\semr{z}$.
First of all we may assume that $ord_z(\tilde\alpha)\klg ord_z(\tilde\beta)$ because $g$ and $g\cdot(\begin{smallmatrix}
0&1\\1&0
\end{smallmatrix})$ represent the same coset in $\Lambda(R/I)$.
Let $m:=ord_z(\tilde\alpha)$ and $n:=\min(ord_z(\tilde\gamma), ord_z(\tilde\delta))$, then we set $h:=(\begin{smallmatrix}z^{-m}&0\\ 0&z^{-n} \end{smallmatrix})$ and we claim that $\overline{hg}=\overline{id}$ in $\Lambda(R/I)$, which means $hg\in Gl_2(R/I\sem{z})$.\\
So we set $hg=
\left(\begin{smallmatrix}
z^{-m}\tilde\alpha&z^{-m}\beta\\
z^{-n}\tilde\gamma &z^{-n}\tilde\delta
\end{smallmatrix}\right)
=:\left(\begin{smallmatrix}
\alpha&\beta\\
\gamma&\delta
\end{smallmatrix}\right)
$ as well as $\tau_\m:=(hg)^{-1}b\sigma^{*}(hg)$ and compute:
\begin{align*}
\tau_M &=\frac{1}{\alpha\delta-\gamma\beta}
    \begin{pmatrix}
    \delta&-\beta\\ -\gamma &\alpha
    \end{pmatrix}
    \begin{pmatrix}
        z & 0 \\ 0 & 1
    \end{pmatrix}
    \begin{pmatrix}
    \leftidx{\sigma}{(\alpha)}&\leftidx{\sigma}{(\beta)}\\
    \leftidx{\sigma}{(\gamma)}&\leftidx{\sigma}{(\delta)}
    \end{pmatrix}\\
 &   =\frac{1}{\alpha\delta-\gamma\beta}
    \begin{pmatrix}
        z\delta\leftidx{\sigma}{(\alpha)}-\beta\leftidx{\sigma}{(\gamma)}
        &z\leftidx{\sigma}{(\beta)}\delta-\leftidx{\sigma}{(\delta)}\beta\\
        -z\gamma\leftidx{\sigma}{(\alpha)}+\alpha\leftidx{\sigma}{(\gamma)}
        &-z\leftidx{\sigma}{(\beta)}\gamma+\leftidx{\sigma}{(\delta)}\alpha
    \end{pmatrix}
\end{align*}

The choice of $m$ and $n$ directly implies $0=ord_z(\alpha)\klg ord_z(\beta), ord_z(\gamma), ord_z(\delta)$ and in addition $ord_z(\gamma)=0$ or $ord_z(\delta)=0$. Consequently, $hg\in Mat_2(R/I\sem{z})$ and it further follows that $ord_z(-z\gamma\leftidx{\sigma}{(\alpha)}+\alpha\leftidx{\sigma}{(\gamma)})=0$ or $ord_z(-z\leftidx{\sigma}{(\beta)}\gamma+\leftidx{\sigma}{(\delta)}\alpha)=0$. By assumption $hg$ is in $\Lambda(R/I)$ and therefore $\tau_M$ has coefficients in $R/I\sem{z}$. But in both cases this is only possible if $ord_z(\alpha\delta-\gamma\beta)\klg 0$ and hence equals 0. So we have $det(hg)\in R/I\sem{z}^{{*}}$ and thus $hg\in Gl_2(R/I\lsem z\rsem )$. \\
This proves in particular that $\Lambda(\F_q)={J}/Stab(id)\simeq\Lambda(R/I)$. Now consider the case that $R/I$ is not an integral domain. 
Since $R/I$ is noetherian, there is only a finite number of minimal prime ideals in $R/I$. Let $\{p_j\}_{j\in P}$ be the set of these minimal prime ideals in $R/I$ and set $A_j:=\bruch{R/I}{p_j}$. As we have remarked the maps $\F_q\xrightarrow{\iota} R/I\xrightarrow{pr_j}A_j$ induce maps:
\begin{equation}\label{aiaj}
\xymatrix{
\Lambda(\F_q)\ar[r]^{\Lambda(\iota)} \ar[dr]_{f_j}^{\widesim}&\Lambda(R/I)\ar[d]^{\Lambda(pr_j)}\\
&\Lambda{(A_j)}
}\qquad\mbox{here $f_j:=\Lambda(pr_j\circ\iota)$}
\end{equation}
The group ${J}$ operates equivariantly under these maps and from the above comment it follows that the maps $f_j$ are bijective. Therefore $\Lambda(\iota)$ is injective and we will prove that it is also surjective. \\
We claim that $f^{-1}_j\circ\Lambda(pr_j)=f^{-1}_k\circ\Lambda(pr_k)$ for all $j,k\in P$. If $P$ has one element we are done, otherwise we start with some $j_1\in P$ and denote by $\overline{\{p_{j_1}\}}\subseteq Spec\ R/I$ the Zariski closure of the one point set corresponding to the minimal prime ideal $p_{j_1}$. Since $Spec\ R/I$ is connected we find $j_2\in P$ with $j_1\neq j_2$ and $\overline{\{p_{j_1}\}}\bigcap \overline{\{p_{j_2}\}}\neq\emptyset$. Therefore we find a prime ideal $p$ in $R/I$ with $p_{j_1}\subseteq p$ and $p_{j_2}\subseteq p$. The projections of $A_{j_1}$ and $A_{j_2}$ onto $A:=\bruch{R/I}{p}$ and the isomorphism $\Lambda(\F_q)\isom \Lambda(A)$ induced by the injection $\F_q\hookrightarrow A$ yield the following diagramm.
\begin{equation*}
\xymatrix{
&\Lambda(R/I)\ar[ld]^{\Lambda(pr_{j_1})} \ar[rd]_{\Lambda(pr_{j_2})}& \\
\Lambda(A_{j_1}) \ar[rdd]_{f_{j_1}^{-1}} \ar[rd] &&\Lambda(A_{j_2}) \ar[ldd]^{f_{j_2}^{-1}} \ar[ld] \\
&\Lambda(A)\ar[d]&\\
&\Lambda(\F_q)&
}
\end{equation*}
It is clear that the lower triangles and the upper square commute, which proves the assertion for $j_1\neq j_2$. Next we choose $j_3\in P$ such that $\overline{\{p_{j_1}\}}\bigcap \overline{\{p_{j_3}\}}\neq\emptyset$ or $\overline{\{p_{j_2}\}}\bigcap \overline{\{p_{j_3}\}}\neq\emptyset$. The same argument proves the assertion for $j\neq k$ with $j,k\in\{j_1,j_2,j_3\}$. Since $Spec\ R/I$ is connected we can continue in the same way until the claim is proved for all $j,k\in P$.\\
The diagram \eqref{aiaj} implies that the (now well-defined) map $f:=f_{j}^{-1}\circ\Lambda(pr_j):\Lambda(R/I)\to \Lambda(\F_q)$ is surjective, we prove now that it is injective. Since ${J}$ operates transitively on $\Lambda(\F_q)$ and equivariantly under the map $f$ it suffices to prove that the preimage of $\overline{id}\in\Lambda(\F_q)$ under the map $f$ equals $\{\ \overline{id}\ \}\in \Lambda(R/I)$. So choose $g\in\Lambda(R/I)$ with $f(g)=\overline{id}$.
We have $f_{j}(\overline{id})=\overline{id}\in \Lambda(A_j)$, which implies together with the above commutative diagram $g\ \mod p_j\in Gl_2(A_j\sem{z})$ for all $j\in P$. Since $\ds\bigcap_{j\in P}p_j=\sqrt{(0)}=(0)$ in $R/I$, we see that $g\in Mat_2(R/I\sem{z})$ and since the same is true for $g^{-1}$ we see $g\in Gl_2(R/I\sem{z})$, which proves $\overline{g}=\overline{id}\in \Lambda(R/I)$ and therefore the injectivity of $f$. We conlude that $\Lambda(\iota)$ is bijective, which ends the proof of the lemma.
}\\

Now this lemma and the bijection \eqref{bilam} enable us to determine the $R/I$-valued points of $\gf$. We have already computed
\begin{equation*}
{J}
=\left\{h\in Gl_2(\F_q^{alg}\semr{z}):hb=b\sigma^{*} h\right\}=\begin{pmatrix}
\F_q\semr{z}^{*} &0\\ 0&\F_q\semr{z}^{*}
\end{pmatrix}.
\end{equation*}
The stabilizer of $id\in\Lambda(R/I)$ is given by 
\begin{equation*}
Stab(id)=\{h\in {J}:\overline{h\cdot id}=\overline{id}\ in\ \Lambda{(R/I)}\}=
    \begin{pmatrix}
        \F_q\lsem z\rsem ^{*} & 0\\ 0 & \F_q\lsem z\rsem ^{*}
    \end{pmatrix}.
\end{equation*}
Therefore,
$\left\{\big(
\begin{smallmatrix}
z^i&0\\0&z^j
\end{smallmatrix}\big)|\ i,j\in \Z
\right\}$ is a system of representatives of ${J}/Stab(id)$ and we get the following bijection
\begin{equation}\label{kwertigebijektion}
\bruch{{J}}{Stab(id)}\to \gf(Spec\ R/I),\quad
(\begin{smallmatrix}
z^i & 0 \\ 0 & z^j
\end{smallmatrix}
)\mapsto 
\big(R/I\lsem z\rsem ^2,
(\begin{smallmatrix}
z & 0 \\ 0 & 1
\end{smallmatrix}),\ 
(\begin{smallmatrix}
z^i & 0 \\ 0 & z^j
\end{smallmatrix})
\big).
\end{equation}

In particular for every arbitrary field $k$ over $\F_q$ it holds true that every $k$-valued point of $\gf$ is already an $\F_q$-valued point. This implies that a formal scheme that ind-represents $\gf$ must be $0$-dimensional and therefore of the form $\coprod_{(i,j)\in\Z^2}\ Spf\left(A_{i,j}\right)$, with $Spf\ A_{i,j}$ consisting of one point and having $Spec\ \F_q$ as reduced subscheme. For $R$ connected, $I:=\sqrt{(0)}$ in $R$ as before and $pr:R\to R/I$ we want to define bijective maps $nat(R)$ as in the diagram.

\begin{equation}\label{achterdiagramm}
\xymatrix{
\ds\coprod_{(i,j)\in \Z^2} Hom(Spec\ R,\ Spf \ A_{i,j}) \ar[d]^{\coprod Hom(pr,\ Spf \ A_{i,j})}
\ar[rr]^{\quad nat(R)}
&& \gf(Spec\ R) \ar[d]^{\gf(pr)}\\
\ds\coprod_{(i,j)\in \Z^2} Hom(Spec\ R/I,\ Spf \ A_{i,j}) 
%
\ar[rr]^{\qquad nat(R/I)}
&& \gf(Spec\ R/I)\\
(f:Spec\ R/I\to Spf\ A_{ij})\ar@{|->}[rr]
&& 
x_{ij}
}\end{equation}

Here $x_{ij}$ denotes the point $\Big(R/I\sem{z}^2,({\begin{smallmatrix}
z & 0\\ 0 & 1
\end{smallmatrix}}),({\begin{smallmatrix}
z^i & 0\\ 0 & z^j
\end{smallmatrix}})
\Big)$ in $\gf(Spec\ R/I)$. From the lower map we know that it is bijective, since $Hom(Spec\ R/I,\ Spf \ A_{i,j})$ consists of exactly one element. This is because a morphism in this set factors uniquely through the reduced subscheme, which means $Hom(Spec\ R/I,\ Spf \ A_{i,j})=Hom(Spec\ R/I, Spec\ \F_q)$.

We set $pr:R\to R/I$ and define:
\nomenclature{\sortkey{Wr}$W_{i,j}(R)$}{is defined by $W_{i,j}(R):=\left\{x\in \gf(Spec\ R)\ \vert\ \gf(pr)(x)=x_{ij}\right\}$}
\begin{equation*}
W_{i,j}(R):=\left\{x\in \gf(Spec\ R)\ |\ \gf(pr)(x)=x_{ij}\right\}
\end{equation*}
Since we have $\gf(Spec\ R/I)=\{x_{ij}\ |\ (i,j)\in\Z^2\}$, the disjoint union $\coprod_{(i,j)\in Z^2} W_{i,j}(R)$ is equal to $\gf(Spec\ R)$.
Now ${J}$ operates transitively on $\gf(Spec\ R/I)$ and equivariantly under the maps $\gf(pr)$. 
Therefore, we see that ${J}\times W_{0,0}(R)=\gf(Spec\ R)$ or even more precisely $\left\{\big(\begin{smallmatrix}z^n & 0\\ 0& z^m\end{smallmatrix}\big)\right\}\times W_{i,j}(R)=W_{i+n,j+m}(R)$. The maps $nat(R)$ have to make the above diagram commutative. Since $\gf(pr)(W_{i,j}(R))=\{x_{ij}\}$ and since $x_{ij}$ is also the image of the unique morphism $f:Spec\ R/I\to Spf\ A_{ij}$, this is only possible if $nat(R)$ maps $Hom(Spec\ R,Spf\ A_{ij})$ bijectively to $W_{i,j}(R)$.

The operation \eqref{joperation}  of ${J}$ and the bijection \eqref{bilam} show us that the isomorphism
 $W_{0,0}(R)\isom W_{i,j}(R)$ is given by
 \begin{equation*}
 \Big(M,\tau_M,\eta\Big)\mapsto \Big(M,\tau_M,\big(\begin{smallmatrix}
z^i & 0 \\ 0 & z^j
\end{smallmatrix}\big) \cdot \eta\Big)\quad .
 \end{equation*}
In particular, it follows that $Hom(Spec\ R,\ Spf\ A_{0,0})\simeq Hom(Spec\ R,\ Spf\ A_{i,j})$ for all $i,j\in \Z$, which means that all the adic rings $A_{i,j}$ are isomorphic.
Before we are able to define the maps $nat(R)$ we have to determine the sets $\gf(Spec\ R)$. In fact we will now inductively determine the sets $\gf(Spec\ R/I^{q^n})$ and
because of the above isomorphism, it suffices to determine $W_{0,0}(R/I^{q^n})$.
Let $h\in I$, then we define $x_n(h)$ to be the tripel:
\begin{equation}\label{xnh}
x_n(h):=\left(M_n=
R/I^{q^n}\lsem z\rsem ^2, 
\tau_{M_n}=\begin{pmatrix}
z-\zeta & 0 \\ h & 1
\end{pmatrix},\ 
\eta_{n,h}= 
\begin{pmatrix}
\prod_{i=0}^{n-1}\frac{z}{z-\zeta^{q^i}} & 0 \\
-\sum_{i=0}^{n-1}\frac{h^{q^i}}{(z-\zeta)\dots(z-\zeta^{q^i})} & 1
\end{pmatrix}\right)
\end{equation} 

We verify that $x_n(h)$ represents a point in $W_{0,0}(R/I^{q^n})$. 
First of all the local shtuka is bounded by $\mu=(1,0)$, since $\tau_{M_n}$ has coefficients in $R/I^{q^n}\lsem z\rsem $ and $\det\tau_{M_n}\in (z-\zeta)R/I^{q^n}\lsem z\rsem ^{*}$. 
\nomenclature{\sortkey{xnh}$x_n(h)$}{a special point in $W_{0,0}(R/I^{q^n})$}
The morphism $\eta_{n,h}$ defines a quasi-isogeny $\underline M\to \underline \M_{R/I^{q^n}}$ since 
\begin{align}
\eta_{n,h}^{-1}\ \tau_{\M}\ \sigma^{*}\eta_{n,h}=&
\begin{pmatrix}
\prod_{i=0}^{n-1}\frac{z-\zeta^{q^i}}{z} & 0 \\
\sum_{i=0}^{n-1}\frac{h^{q^i}(z-\zeta^{q^{i+1}})\dots(z-\zeta^{q^{n-1}})}{z^n} & 1
\end{pmatrix}
\begin{pmatrix}
z&0 \\ 0 & 1
\end{pmatrix}
\begin{pmatrix}
\prod_{i=0}^{n-1}\frac{z}{z-\zeta^{q^{i+1}}} & 0 \\
-\sum_{i=0}^{n-1}\frac{h^{q^{i+1}}}{(z-\zeta^q)\dots(z-\zeta^{q^{i+1}})} & 1
\end{pmatrix}\notag\\
=&\begin{pmatrix}
z-\zeta & 0 \\ h & 1
\end{pmatrix}\quad .\label{isshtuka}
\end{align}
The reduction $mod$ $I$ gives $\gf(pr)(x_n(h))=\Big(
R/I\lsem z\rsem ^2, \tau_{M_1}=(\begin{smallmatrix}
z & 0 \\ 0 & 1
\end{smallmatrix}),\ 
\eta_{1,0}= (
\begin{smallmatrix}
1 & 0 \\
0 & 1
\end{smallmatrix})\Big)$. Together this proves that $x_n(h)\in W_{0,0}(R/I^{q^n})$. Furthermore, we see that for $h\neq \tilde h$ the points $x_n(h)$ and $x_n(\tilde h)$ are not isomorphic. This is because
\begin{equation*}
\eta_{n,\tilde h}^{-1}\cdot \eta_{n,h}= 
\begin{pmatrix}
\prod_{i=0}^{n-1}\frac{z-\zeta^{q^i}}{z} & 0 \\
(\prod_{i=0}^{n-1}\frac{z-\zeta^{q^i}}{z})
\sum_{i=0}^{n-1}\frac{\tilde{h}^{q^i}}{(z-\zeta)\dots(z-\zeta^{q^i})}& 1
\end{pmatrix}\cdot\eta_{n,h}=
\begin{pmatrix}
1 & 0 \\
\sum_{i=0}^{n-1}\frac{(\tilde{h}-h)^{q^i}}{(z-\zeta)\dots(z-\zeta^{q^i})}& 1
\end{pmatrix}
\end{equation*}
is not an element in $Gl_2(R/I^{q^n}\sem{z})$.
We will now prove the following crucial lemma inductively.
\lem{}{
Let $x\in W_{0,0}(R/I^{q^n})
$ then $x$ is isomorphic to a point $x_n(h)$ for some $h\in I$.
}

{\it Proof:} The base clause $n=0$ is clear since the set $W_{0,0}(R/I)$ consists of only one point by the previous diagram. 

\textsl{Induction step $n\to n+1$:}
Let $(M,\tau_M,\eta)$ be an arbitrary point in $W_{0,0} 
(R/I^{q^{n+1}})$.
With the bijection \eqref{bilam} this corresponds exactly to the element $\eta\in\Lambda(R/I^{q^{n+1}})$ and it is sufficient to show that there exists a $g\in Gl_2(R\sem{z})$ with $\eta\equiv \eta_{n,\tilde h}\circ g\ \mod I^{q^{n+1}}$ or equally $g\eta^{-1}\equiv\eta_{n,\tilde h}^{-1}
\ \mod I^{q^{n+1}}$ for some $\tilde h\in I$.
By the induction hypothesis we can assume $\eta^{-1}\equiv
\begin{pmatrix}
\prod_{i=0}^{n-1}\left(\frac{z-\zeta^{q^i}}{z}\right) & 0 \\
\sum_{i=0}^{n-1}\frac{h^{q^i}(z-\zeta^{q^{i+1}}) \dots (z-\zeta^{q^{n-1}})}{z^n} & 1
\end{pmatrix}\ \mod I^{q^n}$ and we write
$\eta^{-1}=
\begin{pmatrix}
\prod_{i=0}^{n-1}\left(\frac{z-\zeta^{q^i}}{z}\right) +a& b \\
\sum_{i=0}^{n-1}\frac{h^{q^i}(z-\zeta^{q^{i+1}}) \dots (z-\zeta^{q^{n-1}})}{z^n}+c & 1+d
\end{pmatrix} \ \mod I^{q^{n+1}}$ with representatives $a,b,c,d\in I^{q^n}R\semr{z}$. Then we have:
\begin{align*}
 &\tau_M:=
 \begin{pmatrix}
 \tau_{11}& \tau_{12}\\
 \tau_{21}& \tau_{22}
 \end{pmatrix}
 :=\ \eta^{-1}\tau_\M\sigma^{*}\eta \\
& \equiv 
\begin{pmatrix}
\prod_{i=0}^{n-1}\left(\frac{z-\zeta^{q^i}}{z}\right) +a& b \\
\sum_{i=0}^{n-1}\frac{h^{q^i}(z-\zeta^{q^{i+1}}) \dots (z-\zeta^{q^{n-1}})}{z^n}+c & 1+d
\end{pmatrix}
\begin{pmatrix}
z&0\\ 0 & 1
\end{pmatrix}
\begin{pmatrix}
\prod_{i=0}^{n-1}\left(\frac{z}{z-\zeta^{q^{i+1}}}\right) & 0 \\
-\sum_{i=0}^{n-1}\frac{h^{q^{i+1}}}{(z-\zeta^{q}) \dots (z-\zeta^{q^{i+1}})} & 1
\end{pmatrix}\\
&\equiv  
\begin{pmatrix}
\frac{z(z-\zeta)}{z-\zeta^{q^{n}}}+az\prod_{i=0}^{n-1}\left(\frac{z}{z-\zeta^{q^{i+1}}}\right)-b\sum_{i=0}^{n-1}\frac{h^{q^{i+1}}}{(z-\zeta^{q}) \dots (z-\zeta^{q^{i+1}})}
%
& b \\
z\sum_{i=0}^{n-1}\frac{h^{q^i}(z-\zeta^{q^{i+1}}) \dots (z-\zeta^{q^{n-1}})}{(z-\zeta^{q}) \dots (z-\zeta^{q^{n}})}+cz\prod_{i=0}^{n-1}\left(\frac{z}{z-\zeta^{q^{i+1}}}\right)-(1+d)\sum_{i=0}^{n-1}\frac{h^{q^{i+1}}}{(z-\zeta^{q}) \dots (z-\zeta^{q^{i+1}})} 
& 1+d
\end{pmatrix}\\ & \mod I^{q^{n+1}}
\end{align*}
The equation \eqref{isshtuka} shows $\tau_M\equiv(\begin{smallmatrix}
z-\zeta &0 \\ h & 1
\end{smallmatrix})\mod I^{q^n}
$.
By assumption, the local shtuka $(M,\ \tau_M)$ is bounded by $(1,0)$. In particular, we have $b, 1+d\ \mod I^{q^{n+1}}\in R/I^{q^{n+1}}\sem{z}$. 
Since we want to find a matrix $g$ with $g\eta^{-1}\equiv \eta_{n,h}^{-1}\ \mod I^{q^{n+1}}$ we can assume without loss of generality that the lifts $b,1+d$ lie in $R\sem{z
}$.
We can therefore write $1+d=1+\sum_{i=0}^\infty\delta_iz^i$ with $\delta_i\in I^{q^n}$ for $i\grg 0$. Replacing $z$ by $z-\zeta+\zeta$ shows that we can also write $1+d=1+a_0+ (z-\zeta)a_1$ with $a_0\in I^{q^n}$ and $a_1\in I^{q^n}R\sem{z}$. 
Since $(1+a_0)^{-1}=(1-a_0+a_0^2-a_0^3+\dots)$ is invertible in $R$, we conclude that $1+d$ is invertible in $R\sem{z}$. Concretely we have $(1+d)^{-1}=(1+a_0)^{-1}
\left(1+\sum_{i=0}^{\infty}(-1)^i\big(\frac{(z-\zeta)a_1}{1+a_0}\big)^i\right)$.
\\
Next we know also by the boundedness condition that $\det\tau_M\in (z-\zeta) R/I^{q^{n+1}}\sem{z}^{*}$. Thus, we can find a unit $u\in R\sem{z}^{*}$ with $(u\mod I^{q^{n+1}}) = \frac{\det\tau_M}{z-\zeta}=\frac{\det\eta^{-1}\cdot z\cdot\det\sigma^{*}\eta}{z-\zeta}$ and it follows that
\begin{equation}\label{ukong}
u^{-1}\det\eta^{-1}\equiv\frac{z-\zeta}{z\cdot \det\sigma^{*}\eta} \equiv\prod_{i=0}^{n}\frac{z-\zeta^{q^i}}{z}\mod I^{q^{n+1}} \quad .
\end{equation}
Now we define $g_1\in Gl_2(R\sem{z})$ to be the matrix on the left in the following equation and compute:
\begin{align*}
&\underbrace{
\begin{pmatrix}
(1+d)u^{-1} & -bu^{-1}\\
0 & (1+d)^{-1}
\end{pmatrix}
}_{=:g_1\in Gl_2(R\lsem z\rsem )}
\underbrace{\begin{pmatrix}
\prod_{i=0}^{n-1}\left(\frac{z-\zeta^{q^i}}{z}\right) +a& b \\
\sum_{i=0}^{n-1}\frac{h^{q^i}(z-\zeta^{q^{i+1}}) \dots (z-\zeta^{q^{n-1}})}{z^n}+c & 1+d
\end{pmatrix}}_{\equiv\ \eta^{-1}\mod I^{q^{n+1}}}
\\
=&
\begin{pmatrix}
u^{-1}\det\eta^{-1} & 0 \\
(1+d)^{-1}\left(\sum_{i=0}^{n-1}\frac{h^{q^i}(z-\zeta^{q^{i+1}}) \dots (z-\zeta^{q^{n-1}})}{z^n}+c \right)&1
\end{pmatrix}
\end{align*}
We write $(1+d)^{-1}=d_0+(z-\zeta)d_1$ with $d_1\in 
R\lsem z\rsem $ and $d_0\in R$. By the previous computation of $(1+d)^{-1}$ we know that $d_0\equiv (1+a_0)^{-1}\equiv 1\mod I^{q^n}$ and $d_1\equiv 0\mod I^{q^n}$.
Since $\tau_M 
$ has coefficients in $R/I^{q^{n+1}}\lsem z\rsem $ 
and $\tau_{21}\equiv h\mod I^{q^n}$ by equation \eqref{isshtuka},
we can choose furthermore $w_0\in R$ and $w_1\in R\lsem z\rsem $ with $\tau_{21}\equiv w_0+(z-\zeta)w_1\mod I^{q^{n+1}}$ and $w_0\equiv h\mod I^{q^n}$ as well as $w_1\equiv 0 \mod I^{q^n}$. This is done in the same way as we chose $a_0$ and $a_1$.
We now define the desired $\tilde h\in I$ as:
\begin{equation}
\tilde{h}:=w_0d_0\quad\mbox{\ \ $\Rightarrow$ }\qquad\tilde{h}\equiv h\mod I^{q^n}\quad \Rightarrow\quad \tilde{h}^q\equiv h^q\mod I^{q^{n+1}}\label{hdef}
\end{equation}

The following equivalences then end with the definition of $\widetilde{\eta_{21}}$:

\begin{align}
&\tau_{21}\equiv w_0+(z-\zeta)w_1 \qquad \mod I^{q^{n+1}}
\mbox{\hspace{5.5cm} }\Bigg|\frac{1}{z^{n+1}}\prod_{i=0}^{n-1}z-\zeta^{q^{i+1}} \notag\\
\Leftrightarrow\ & c+\sum_{i=0}^{n-1}\frac{h^{q^i}(z-\zeta^{q^{i+1}}) \dots (z-\zeta^{q^{n-1}})}{z^n}-\left(\prod_{i=0}^{n-1}\frac{z-\zeta^{q^{i+1}}}{z}\right)\frac{(z-\zeta)w_1}{z}\notag\\
& \qquad\quad
\equiv
\frac{w_0}{z^{n+1}}\prod_{i=0}^{n-1}(z-\zeta^{q^{i+1}})+(1+d)\sum_{i=0}^{n-1}\frac{h^{q^{i+1}}(z-\zeta^{q^{i+2}}) \dots (z-\zeta^{q^n})}{z^{n+1}}
\quad \mod I^{q^{n+1}}
 \notag \\
\Leftrightarrow\ & (1+d)^{-1}\left(c+\sum_{i=0}^{n-1}\frac{h^{q^i}(z-\zeta^{q^{i+1}}) \dots (z-\zeta^{q^{n-1}})}{z^n}
-w_1\prod_{i=0}^n\frac{z-\zeta^{q^i}}{z}
\right) 
-\frac{(z-\zeta)d_1w_0}{z^{n+1}}\prod_{i=0}^{n-1}(z-\zeta^{q^{i+1}})\notag \\&
\qquad\quad \equiv\sum_{i=1}^n \frac{h^{q^i}(z-\zeta^{q^{i+1}}) \dots (z-\zeta^{q^{n}})}{z^{n+1}}+\frac{d_0w_0}{z^{n+1}}\prod_{i=0}^{n-1}(z-\zeta^{q^{i+1}}) 
\qquad \quad \mod I^{q^{n+1}}
\notag \\
\underset{\eqref{hdef}}{\Longrightarrow} & 
\underbrace{(1+d)^{-1}\left(c+\sum_{i=0}^{n-1}\frac{h^{q^i}(z-\zeta^{q^{i+1}})\dots(z-\zeta^{q^{n-1}})}
{z^n}+((1+d)w_0d_1-w_1)\prod_{i=0}^n\frac{z-\zeta^{q^i}}{z}
\right)}_{=:\widetilde{\eta_{21}}} \notag \\
&\qquad\quad \equiv \sum_{i=0}^n\frac{\tilde{h}^{q^i}(z-\zeta^{q^{i+1}})\dots(z-\zeta^{q^n})}{z^{n+1}} 
\qquad \mod I^{q^{n+1}} \label{monster}
\end{align}
We define a second element $g_2\in Gl_2(R\sem{z})$ as the left matrix in the following equation and do a final computation for the proof of the lemma:

\begin{align*}
g_2g_1\eta^{-1}\equiv &
\underbrace{
\begin{pmatrix}
1 & 0 \\ w_0d_1-(1+d)^{-1}w_1 & 1
\end{pmatrix}}_{=:g_2\in Gl_2(R\lsem z\rsem )}
\underbrace{
\begin{pmatrix}
u^{-1}\det\eta^{-1} & 0 \\
(1+d)^{-1}\left(\sum_{i=0}^{n-1}\frac{h^{q^i}(z-\zeta^{q^{i+1}})\dots(z-\zeta^{q^{n-1}})}{z^n}+c \right)&1
\end{pmatrix}}_{\equiv g_1\eta^{-1}\ \mod I^{q^{n+1}}}\\
\overset{\eqref{ukong}}{\equiv}&
\begin{pmatrix}
\prod_{i=0}^n\frac{z-\zeta^{q^{i}}}{z} & 0 \\
\widetilde{\eta_{21}}&1
\end{pmatrix}
\overset{\eqref{monster}}{\equiv}
\begin{pmatrix}
\prod_{i=0}^n\frac{z-\zeta^{q^{i}}}{z}& 0\\
\sum_{i=0}^n\frac{\tilde{h}^{q^i}(z-\zeta^{q^{i+1}})\dots(z-\zeta^{q^{n}})}{z^{n+1}} & 1
\end{pmatrix}\quad \mod I^{q^{n+1}}
\end{align*} 

This shows that $g:=g_2g_1$ is the desired matrix since $\eta g_1^{-1}g_2^{-1}\equiv\eta_{n+1,\tilde{h}}\ \mod I^{q^{n+1}}$ implies that $(M, \tau_M, \eta)$ is isomorphic to $x_{n+1}\big(\tilde{h}\big)$, which ends the proof of the lemma.\hfill $\Box$\\

Due to this lemma we now conclude that $\{x_n(h)\ |\ h\in I{/I^{q^n} }\}$ defines a system of representatives of $W_{0,0}(R/I^{q^n})$.

Therefore, we get the following bijection:
\begin{equation}\label{defnatr}
Hom_{\F_q\sem{\zeta}}(\F_q\lsem \zeta, t\rsem ,\ R/I^{q^n}) \longrightarrow 
 W_{0,0}(R/I^{q^n}), \qquad \big(f: t\mapsto h\in I\big)\longmapsto 
x_n(h)
\end{equation}

Since $Spec\ R$ is a noetherian scheme in $\nilp$ , we have $I^m=(0)$ for some $m>>0$. So for $q^n>m$ this means $R/I^{q^n}\simeq R$ and the equation \eqref{defnatr} can be written as

\begin{equation*}
\xymatrix{
Hom(Spec\ R,\ Spf\ \F_q\lsem \zeta, t\rsem )=Hom(Spec\ R/I^{q^n},\ Spf\ \F_q\lsem \zeta, t \rsem )\ar[r] &W_{0,0}(R/I^{q^n}) =W_{0,0}(R)&
}
\end{equation*}

By the operation of ${J}$ this defines also bijections 
$Hom\big(Spec\ R,\ Spf\ \F_q\sem{\zeta, t}\big)\to  W_{i,j}(R)$ and therefore a map 
\begin{equation*}
nat(R):\ Hom(R,\coprod_{i,j\in\Z^2}Spf\ \F_q\lsem\zeta, t\rsem) \to \coprod_{i,j\in\Z^2}W_{i,j}(R)=\gf(Spec\ R)
\end{equation*}
that makes the diagram \eqref{achterdiagramm} commutative. 
These maps define an isomorphism of functors, which shows that the functor $\gf$ is ind-representable by the formal scheme $\gs=\coprod_{i,j\in\Z^2}Spf\ \F_q\lsem\zeta, t\rsem$ and for the proof of the theorem it rests only to compute the universal ind-object.

\subsubsection{The Universal Ind-Object
}\label{universalob}

The functor $\gf$ extends naturally to the category of ind-objects 
by defining $\gf(S):=\varprojlim_n \gf(S_n)$ for an ind-object $S=\varinjlim_n S_n$ with $S_n\in\nilp$. Concretely this means that a point in $\gf(S)$ is given by a family $(M_n, \tau_{M_n}, \eta_n)_{n\in\N}$ with $(M_n, \tau_{M_n}, \eta_n)\in \gf(S_n)$ and $(M_n, \tau_{M_n}, \eta_n)\simeq({M_{n+1}}_{S_n}, \tau_{M_{n+1}}\otimes{id}, \eta_{n+1}\otimes{id})$. 
Now the universal ind-object $x^{univ}$ 
\nomenclature{\sortkey{xuniv}$x^{univ}$}{universal ind-object in $\gf(\gs)$}
is by definition the image of $id_{\gs}$ under the map $nat(\gs):Hom(\gs,\gs)\to \gf(\gs)$.
It has the useful property that it defines the maps $nat(S): Hom(S,\gs)\to \gf(S)$ by $f\mapsto \gf(f)(x^{univ})$.\\
Let $\gs=\coprod_{(i,j)\in\Z^2}U_{ij}$ with $U_{ij}\simeq Spf\ \F_q\sem{\zeta,h}$. It suffices to describe $x^{univ}|_{U_{ij}}$, so we fix $(i,j)\in\Z^2$ and write $A_n=\F_q\sem{\zeta,h}/(\zeta,h)^{q^n}$. Thus, we identify $U_{ij}$ with $\varinjlim_n\ Spec\ A_n$. This identification comes with natural morphisms $\iota_n:Spec\ A_n\to U_{ij}$. They are induced by the projections $\F_q\sem{\zeta,h}\to A_n$.
The definition of $nat(R)$ at the end of the previous section shows that $nat(A_n)$ is given by
\begin{align*}
nat(A_n): \quad &Hom(Spec\ A_n,U_{ij})&\longrightarrow & \quad W_{i,j}(A_n)\\
&\iota_n&\longmapsto & \quad \{\big(\begin{smallmatrix}
  z^i & 0\\ 0 & z^j
  \end{smallmatrix}\big)\}\times x_n(h)\quad .
 \end{align*}
Since we have $$\begin{pmatrix}
z^i&0\\ 0 & z^j
\end{pmatrix}
\begin{pmatrix}
\prod_{i=0}^{n-1}\frac{z}{z-\zeta^{q^i}} & 0 \\
-\sum_{i=0}^{n-1}\frac{{h}^{q^i}}{(z-\zeta)\dots(z-\zeta^{q^i})}& 1
\end{pmatrix}=
\begin{pmatrix}
z^i\prod_{i=0}^{n-1}\frac{z}{z-\zeta^{q^i}} & 0 \\
-z^j\sum_{i=0}^{n-1}\frac{{h}^{q^i}}{(z-\zeta)\dots(z-\zeta^{q^i})}& z^j
\end{pmatrix}$$
the universal ind object $\ds x^{univ}=(x_n^{univ})_{n\in\N}\in \varprojlim_n \gf(\coprod_{i,j\in\Z^2} U_{ij}^n)$ is given by
$$
x_n^{univ}\Big|_{U_{ij}^n}=\left(\F_q\sem{\zeta,h}/(\zeta, h)^{q^n}\sem{z-\zeta}^2,\ \begin{pmatrix}
z-\zeta & 0 \\ h & 1 
\end{pmatrix}
,\ 
\begin{pmatrix}
z^i\prod_{i=0}^{n-1}\frac{z}{z-\zeta^{q^i}} & 0 \\
-z^j\sum_{i=0}^{n-1}\frac{{h}^{q^i}}{(z-\zeta)\dots(z-\zeta^{q^i})}& z^j
\end{pmatrix}\right)\quad .
$$

This proves the theorem. \bo

\subsubsection{Remark on formal schemes locally formally of finite type}
Recall that a formal $\F_q\sem{z}$-scheme is said to be locally formally of finite type if it is locally noetherian and if the induced morphism on the reduced subschemes $S^{red}\to \F_q$ is locally of finite type.  We denote the category of these schemes by $\flfts$. For a formal scheme $S\in\flfts$ we denote by $I_S$ the maximal ideal of definition \cite[10.5.4]{GD71} and by $S_n:=(S,\oo_S/I_S^n)$ the corresponding closed subscheme in $\nilp$ so that we can write $S=\varinjlim S_n$. Since $\flfts$ is a full subcategory of the category of ind-schemes of $\nilp$, we can extend $\gf$ from $\nilp$ to $\flfts$ in the same way as above.
We remark that the points $(M_n, \tau_{M_n}, \eta_n)_{n\in\N}\in\gf(S)$ correspond bijectively to pairs $(M,\tau_M)$ consisting of a locally free $\oo_S\sem{z}$-module $M$ together with an isomorphism $\tau_M:\sigma^{*} M[\frac{1}{z-\zeta}]\to M[\frac{1}{z-\zeta}]$. However, one has to be careful with the quasi-isogenies $\eta_n$ because it is \underline{not} true that they define an isomorphism $\eta:M[\frac{1}{z-\zeta}]\to \M_S[\frac{1}{z-\zeta}]$. This can be seen for example from the universal object that we computed above. We will prove in lemma \ref{etalim} that the assertion becomes true for $(\eta_n)_{n\in\N}$ if we pass to the associated module over the rigid analytic space $S^{rig}$.

\section{$z$-Isocrystals with Hodge-Pink Structures}\label{ziso}
Let $\underline\M$ be a fixed local shtuka of rank $r$ over $\F_q$ and $(\mu_1,\dots,\mu_r)\in\Z^r$ with $\mu_1\grg
\dots\grg\mu_r$ and $\gf$ be the associated Rapoport-Zink functor. 
The main goal of this section is to associate with a point $x\in \gf(S)$ a Hodge-Pink lattice on the $z$-isocrystal $\underline\M_{\F_q\semr{z}}$. 
The functor $rig$ associates with every formal scheme $S\in \flfts$ over $Spf\ \F_q\sem{z}$ its "generic fibre" $S^{rig}$, which is a rigid analytic space over $\F_q\semr{z}$. 
Bosch defines in his book \cite{Bos14} the functor $rig$ only for formal schemes that are locally topoloigically of finite type. Since the Rapoport-Zink space $\gs$ is not locally topologically of finite type, we need a more general construction by Berthelot that defines $rig$ on $\flfts$ and which is described in \cite[5.5]{RZ96}. Note that if an affine formal scheme $S=Spf\ A$ is not topologically of finite type, the generic fibre $S^{rig}$ will not be an affinoid space any more. Nevertheless it is still a quasi-stein space, that is, there exists an open admissible covering by 
affinoid subspaces $U_1,U_2,\dots \subseteq S^{rig}$ such that the image of $\oo_X(U_{i+1})$ is dense in $\oo_X(U_i)$. If $I$ is a defining ideal for $A$ and $f_1,\dots, f_r$ is a system of generators of $I$, then the covering $U_n$ is given by $U_n=Sp(B_n\otimes_{\F_q\sem{\zeta}} {\F_q\semr{\zeta}})$, where $B_n$ are $F_q\sem{z}$ algebras topologically of finte type \cite[5.5]{RZ96} given by
\begin{equation}
\nomenclature{\sortkey{Bn}$B_n$}{Ring used for the construction of the special fiber of a formal scheme}
\label{ringbn}
B_n:=A\grk{T_1,\dots ,T_r}/(f_1^n-\zeta T_1,\dots,f_r^n-\zeta T_r).
\end{equation}
For a rigid analytic space $X$ over $\F_q\semr{\zeta}$ we define similar as before $\oo_X\sem{z-\zeta}$ to be the sheaf of rings for the strong Grothendieck topology on $X$ that associates with every admissible open subset $U\subseteq X$ the ring $\oo_X(U)\sem{z-\zeta}$. Note that in contrast to the sheaves $\oo_S\sem{z}=\oo_S\sem{z-\zeta}$ on some formal scheme $S\in \flfts$, we cannot write $\oo_X\sem{z}$ for the sheaf $\oo_X\sem{z-\zeta}$, because $\zeta$ is invertible in $\oo_X(U)$.
Furthermore, let $\oo_X\semr{z-\zeta}$
\nomenclature{\sortkey{Oxa}$\oo_X\sem{z-\zeta}$}{sheaf on a rigid analytic space $X$ given by $U\mapsto \oo_X(U)\sem{z-\zeta}$}
\nomenclature{\sortkey{Oxb}$\oo_X\semr{z-\zeta}$}{sheaf on a rigid analytic space $X$ associated with the presheaf $U\mapsto \oo_X(U)\sem{z-\zeta}[\frac{1}{z-\zeta}]$}
 be the sheaf associated with the presheaf given by $U\mapsto \oo_X(U)\sem{z-\zeta}[\frac{1}{z-\zeta}]$.

There is a natural map of ringed spaces $sp:S^{rig}\to S$, which induces a morphism of sheaves $sp^{-1}\oo_S\sem{z}\to \oo_{S^{rig}},\quad z\mapsto \zeta$ as well as $sp^{-1}\oo_S\sem{z}\to \oo_{S^{rig}}\sem{z-\zeta},\quad z\mapsto \zeta+(z-\zeta)$ and $sp^{-1}\oo_S\semr{z-\zeta}\to \oo_{S^{rig}}\semr{z-\zeta}$ on $S^{rig}$.
Now let $M$ be an $\oo_S\sem{z}$-module so that $sp^{-1}(M)$ is an $sp^{-1}\oo_S\sem{z}$-module. By abuse of notation we write  $M\otimes_{\oo_S\sem{z}}\oo_{S^{rig}}$ instead of $sp^{-1}(M)\otimes_{sp^{-1}\oo_S\sem{z}}\oo_{S^{rig}}$ and even shorter we write $M_{S^{rig}}\sem{z-\zeta}:=M\otimes_{\oo_S\sem{z}}\oo_{S^{rig}}\sem{z-\zeta}$ and $M_{S^{rig}}\semr{z-\zeta}:=M\otimes_{\oo_S\sem{z}}\oo_{S^{rig}}\semr{z-\zeta}$ for the respective tensor products.

In the following subsection \ref{neuziso} we define $z$-isocrystals, Hodge-Pink lattices and their moduli spaces.  
In subsection \ref{associatedhps} we then describe how one associates with a point $x\in\gf(S)$ a Hodge-Pink lattice over $S$ on ${\underline\M}_{\F_q\semr{z}}$.

\subsection{Definition of $z$-Isocrystals with Hodge-Pink Structures and their Moduli Spaces}
\label{neuziso}

Let $k$ be an arbitrary field over $\F_q$, later we will only be interested in the case $k=\F_q$.

\df{{\cite[Def. 5.1]{HK15}}}{\mbox{}\label{zisodef}
A $z$-isocrystal over $k$ is a pair $(D,\tau_D)$ consisting of a finite dimensional $k\semr{z}$ vector space together with a $k\semr{z}$-isomorphism $\tau_D:\sigma^{*} D\isom D$. \\ A morphism $(D,\tau_D)\to (D',\tau_{D'})$ of $z$-isocrystals is a $k\semr{z}$-homomorphism $f:D\to D'$ satisfying $\tau_{D'}\circ \sigma^{*} f=f\circ \tau_D$.}

The associated $z$-isocrystal of the local shtuka $\underline \M$ over $\F_q$ is given by $\underline\M_{\F_q\semr{z}}:=(\M\otimes_{\F_q\sem{z}}\F_q\semr{z}, \tau_{\M}\otimes id)$.\smallskip

\df{{compare also \cite[D\'efinition 3.9]{GL11}}}{\mbox{}\label{hodgeneu}
Let $X$ be a rigid analytic $k\semr{\zeta}$-space. A Hodge-Pink lattice over $X$ on $(D,\tau_D)$ is a locally free $\oo_X\sem{z-\zeta}$-module $\q$ with an isomorphism
%
$\q\otimes_{\oo_X\sem{z-\zeta}}\oo_X\semr{z-\zeta}\simeq 
D\otimes_{k\semr{z}}\oo_X\semr{z-\zeta}$ satisfying the condition that there exists $d,e\in\Z$ such that 
\begin{equation*}
(z-\zeta)^{-d}\p_X\subseteq \q\subseteq (z-\zeta)^{-e}\p_X
\end{equation*}
where $\p_X=D\otimes_{k\semr{z}}\oo_X\sem{z-\zeta}$. The Hodge-Pink lattice is then said to be of amplitude $[d,e]$.\\
For a formal $k\sem{\zeta}$-scheme $S$, locally formally of finite type, a Hodge-Pink lattice over $S$ on $(D,\tau_D)$ is a Hodge-Pink lattice over $S^{rig}$.}

A tripel $(D,\tau_D,\q)$ is called a $z$-isocrystal with Hodge-Pink lattice over $X$. A morphism $(D,\tau_D,\q)\to (D',\tau_{D'},\q')$ of two $z$-isocrystals with Hodge-Pink lattice over $X$ is a morhpism $f:(D,\tau_D)\to (D',\tau_{D'})$ of $z$-isocrystals such that $(f\otimes id)(\q)\subseteq \q'$.

The module $\p_X$ is itself a Hodge-Pink lattice and is called the special Hodge-Pink lattice over $X$.

Let $X$ be a quasi-stein space and $B=\oo_X(X)$ be its ring of integers. It follows by \cite[Satz 2.4]{Kie66} that there is a fully faithful functor $-\otimes_B\oo_X:\{B-modules\}\to \{\oo_X-modules\}$ and respective functors from the categories of $B\sem{z-\zeta}$- and $B\semr{z-\zeta}$-modules to the categories of $\oo_X\sem{z-\zeta}$- and $\oo_X\semr{z-\zeta}$-modules. Using this fact as well as \cite[Lemma 2.2.3]{Sch14} and \cite[Proposition 2.2.5]{Sch14} one can see that the Hodge-Pink lattices over $B$ correspond bijectively to the Hodge-Pink lattices over $X$.

\df{{compare \cite[2.3]{Sch14}}}{\mbox{}
Let $(D,\tau_D)$ be a $z$-isocrystal of rank $r$ and $\omega_1,\dots, \omega_r\in \Z$ with $\omega_1\grg\dots\grg\omega_r$. Then a Hodge-Pink lattice $\q$ over a $k\semr{\zeta}$ algebra $R$ (resp. a rigid analytic $k\semr{\zeta}$-space) is said to be bounded by $\omega=(\omega_1,\dots,\omega_r)$ if it satisfies
\begin{equation*}
(z-\zeta)^{\omega_{1}+\dots+\omega_i}\cdot\bigwedge\nolimits^{\!i}\p\subseteq \bigwedge\nolimits^{\!i}\q\qquad \mbox{\ for }1\klg i\klg r\mbox{\quad with equality for $i=r$}.
\end{equation*}
}

By \cite[Lemma 2.3.2]{Sch14} this condition is equivalent to 
\begin{equation*}
\bigwedge\nolimits^{\!i}\q\subseteq(z-\zeta)^{\omega_{r-i+1}+\dots+\omega_r}\cdot\bigwedge\nolimits^{\!i}\p\qquad \mbox{\ for }1\klg i\klg r\mbox{\quad with equality for $i=r$}.
\end{equation*}

\subsubsection{Moduli Spaces for Hodge-Pink Structures}
\nomenclature{\sortkey{Qa}$\mathcal{Q},\ \Qf$}{functors that classify Hodge-Pink structures}
\nomenclature{\sortkey{Qs}$\Qs,\ \Qsr$}{representing spaces of $\mathcal{Q}$ and $\Qf$}
We fix a $z$-isocrystal $\underline D=(D,\tau_D)$ of rank $r$ over $k$ and $\omega=(\omega_1,\dots,\omega_r)\in \Z^r$ with $\omega_1\grg\dots\grg\omega_r$ and define the following functor:
\begin{eqnarray*}
    \mathcal{Q}_{D,\skl\omega}:(k\semr{\zeta}-algebras)    \longrightarrow& Set\\
    R\longmapsto&\left\{
        \minibox{Hodge-Pink lattices $\q$ over $R$ on $(D,\tau_D)$,\\ which are bounded by $\omega$}\right\}\\
    (f:R\to R')\longmapsto & \q\mapsto \q\otimes_{R\sem{z-\zeta}}R'\sem{z-\zeta}
\end{eqnarray*} 

\noindent It is a result in the PhD thesis of Tim Schauch that this functor is locally representable by a projective $k\semr{\zeta}$-scheme \cite[Subsection 2.4]{Sch14}.
The corresponding functor 
\begin{equation*}
\widetilde{\mathcal Q}_{D,\skl\omega}:(rigid\ k\semr{\zeta}-spaces)^{op}    \longrightarrow Set
\quad     
X\longmapsto \left\{ \minibox{Hodge-Pink lattices $\q$ over $X$ on $(D,\tau_D)$,\\ which are bounded by $\omega$}\right\}
 \end{equation*}
 is representable by the rigid analytic space $\mathcal Q_{D,\klg\omega}^{rig}$.
Whenever $\underline D$ and $\omega$ are fixed, we will also write $\mathcal{Q}:=\mathcal{Q}_{D,\skl\omega}$ and $\Qf:=\mathcal Q_{D,\klg\omega}^{rig}=\widetilde{\mathcal Q}_{D,\skl\omega}$.

\bsp{}{\label{beispielq}Let $(D,\tau_D)=\big(\F_q\semr{z}^2,(\begin{smallmatrix}
z&0\\0&1
\end{smallmatrix})\big)
$ and $\omega=(0,-1)$. Then for an $\F_q\sem{\zeta}$-algebra
$R$ we have $D\otimes_{\F_q\semr{z}}R\semr{z-\zeta}=R\semr{z-\zeta}^2$ and:
\begin{equation*}
\mathcal{Q}(R)=\left\{\q\subseteq R\semr{z-\zeta}^2\ |\ \q\mbox{ is a Hodge-Pink lattice bounded by } (0,-1)\right\}
\end{equation*}
We write $\p_R=R\sem{z-\zeta}^2=e_1R\sem{z-\zeta}+e_2R\sem{z-\zeta}$ with the standard vectors $e_1=(\begin{smallmatrix}
1\\0
\end{smallmatrix})
$ and $e_2=(\begin{smallmatrix}
0\\1
\end{smallmatrix})
$. The boundedness condition by $\mu$ means for $i=1$ that 
\begin{equation*}
\q\subseteq (z-\zeta)^{-1}\p_R=(z-\zeta)^{-1}R\sem{z-\zeta}^2\mbox{ and for $i=2$ that}
\end{equation*}
\begin{equation*}
\bigwedge_{R\sem{z-\zeta}}^2\q=(z-\zeta)^{-1}\bigwedge_{R\sem{z-\zeta}}^2\p_R=e_1\wedge e_2\ (z-\zeta)^{-1}R\sem{z-\zeta}.
\end{equation*}
Because of the latter equality there are two elements 
$(\begin{smallmatrix}
r_1\\r_2
\end{smallmatrix})
,(\begin{smallmatrix}
s_1\\s_2
\end{smallmatrix})\in\q
$ with $r_1,r_2,s_1,s_2\in (z-\zeta)^{-1}R\sem{z-\zeta}$ and $(z-\zeta)^{-1}\cdot e_1\wedge e_2=(\begin{smallmatrix}
r_1\\r_2
\end{smallmatrix})
\wedge(\begin{smallmatrix}
s_1\\s_2
\end{smallmatrix})=det\begin{pmatrix}
r_1&s_1\\r_2&s_2
\end{pmatrix}\cdot e_1\wedge e_2
$. This implies $e_1=(z-\zeta)s_2\bvecc{r_1}{r_2}-(z-\zeta)r_2\bvecc{s_1}{s_2}\in\q$ and $e_2=(z-\zeta)s_1\bvecc{r_1}{r_2}-(z-\zeta)r_1\bvecc{s_1}{s_2}\in\q$ and, consequently, $\p_R\subseteq \q$.
Furthermore, the equality condition for $i=2$ implies that the inclusions $\p_R\subsetneq \q\subsetneq (z-\zeta)^{-1}\p_R$ are proper because otherwise we would have $\bigwedge_{R\sem{z-\zeta}}^2\q=e_1\wedge e_2\ (z-\zeta)^{-2}R\sem{z-\zeta}$ or $\bigwedge_{R\sem{z-\zeta}}^2\q=e_1\wedge e_2\ R\sem{z-\zeta}$, which is a contradiction to the boundedness condition for $i=2$. \\
Thus, we can uniquely write $\q=(z-\zeta)^{-1}N+R\sem{z-\zeta}^2$ with some $R$-module $0\neq N\subsetneq R^2$. The condition that $\q$ is a Hodge-Pink lattice implies that $N\simeq\q/\p_R\subset (z-\zeta)^{-1}\p_R/\p_R\simeq R^2$ and $(z-\zeta)^{-1}\p_R/\q\simeq R^2/N$ are locally free $R$-modules and since we have $0\neq N\subsetneq R^2$ they are locally free of rank $1$. Therefore, the functor $\mathcal{Q}_{D,\skl\omega}$ is equivalent to the functor
\begin{eqnarray*}
\{\F_q\semr{\zeta}-algebras\}&&\to Set\\
R&&\mapsto \left\{N\subseteq R^2\Bigg|\ \minibox{N is a locally free $R$-module such that\\ $R^2/N$ is a locally free $R$-module of rank $1$}\right\}.
\end{eqnarray*}

\noindent This is the well-known Grassmannian functor $Grass_{1,2}$ which is representable by the projective space $\PP_{\F_q\semr{\zeta}}^1$ (see \cite[Lemma 8.14, 8.5]{GW10}).
It follows that $\Qf_{D,\klg\omega}$ is representable by $\PP_{\F_q\semr{\zeta}}^{1,rig}$. 
}

\subsection{Hodge-Pink Structures associated with Points in $\gf(S)$}\label{associatedhps}
Now we would like to associate with every point $x$ in $\gf(S)$ a Hodge-Pink lattice over $S$ on the $z$-isocrystal $\underline{\M}_{\F_q\semr{z}}$.
A point $x$ is represented by a compatible family $(M_n, \tau_{M_n}, \eta_{n})_{n\in \N}$. We remarked at the end of the second section that this family provides a locally free $\oo_S\sem{z}$-module M and an isomorphism $\tau_M:\sigma^{*} M[\frac{1}{z-\zeta}]\to M[\frac{1}{z-\zeta}]$.

We saw that the family $(\eta_n)$ of quasi-isognies does not generally define a morphism $\eta:M[\frac{1}{z-\zeta}]\to \M_S[\frac{1}{z-\zeta}]$, but we will need and prove in proposition \ref{etainteger}, that this becomes true after passing to the generic fibre. This means, we have a rigid analytic isomorphism
\begin{equation*}
\displaystyle\sigma^{*}\eta_M:\sigma^{*}M_{S^{rig}}\sem{z-\zeta}\to \sigma^{*}\M_{S^{rig}}\sem{z-\zeta}.
\end{equation*}

Before proving this proposition we need a lemma and the definition of the ring $A\lsem z,z^{-1}\}$ and its Element $l_{-}$.

\subsubsection{The Ring $A\lsem z,z^{-1}\}$ and its Element $l_{-}$}
\label{bizzarring}
Let $S=Spf\ A\in \flfts$ such that $\zeta$ is not nilpotent in $A$ and denote by $I:=I_S(S)$ the maximal defining ideal of the $\F_q\sem{z}$-algebra A. The proof of the next lemma and the next proposition need the notion of the ring $A\lsem z,z^{-1}\}$ and its element $l_{-}$. 
For an element $0\neq a\in A$ we set 
\begin{equation*}
v(a):=\max\{n\in \N_0\ |\ a\in I^n\}\mbox{ and }v(0)=\infty\quad .
\end{equation*}
We have $v(a+b)\grg \min\{v(a),v(b)\}$ and $v(ab)\grg v(a)+v(b)$ and for a sequence $(a_n)_{n\in\N}$ in $A$ the condition $\lim_{n\to \infty}a_n=0$ is equivalent to $\lim_{n\to \infty}v(a_n)=\infty$.
Now we define
\nomenclature{\sortkey{A}$A\semm$}{Ring of Laurent series converging on the open unit disk}
\begin{equation}
A\lsem z,z^{-1}\}:=\left\{\sum_{i=-\infty}^{\infty}a_iz^{i}\ |\ a_i\in A \mbox{ and } \lim_{i\to -\infty} v(a_i)+ni=\infty\quad \forall n>0\right\}.
\end{equation}
Since $\F_q\sem{\zeta}\hookrightarrow A$ we have an inclusion $\F_q\sem{\zeta}\semm\subseteq A\semm$. There is a particular element $l_{-}\in \fzet\semm$ which is defined by 
\begin{equation*}
\nomenclature{\sortkey{l-}$\lper$}{Special element in $\fzet\semm$}
\lper:=\prod_{i\in \N_0}\left(1-\frac{\zeta^{q^i}}{z}\right).
\end{equation*}
This is indeed an element in $\fzet\semm$ since the coefficient $a_{-k}$ of $z^{-k}$ is given by $\ds\sum_{0\klg i_1<\dots<i_k}(-\zeta^{q^{i_1}})\dots(-\zeta^{q^{i_k}})$ with $v(a_{-k})=1+q+\dots +
q^{k-1}=\frac{q^k-1}{q-1}$. We define the rings $B_n$ and the affinoid spaces $U_n$ with respect to $A$ as defined on page \pageref{ringbn}. We have $\oo_{S^{rig}}(U_n)=B_n\otimes_{\F_q\sem{z}}\F_q\semr{z}$ and $S^{rig}=\cup_{n\in\N}U_n$. Since for the open immersions $U_n\hookrightarrow U_{n+1}$ the maps 
$B_{n+1}\otimes_{\F_q\sem{z}}\F_q\semr{z}\to B_n\otimes_{\F_q\sem{z}}\F_q\semr{z},\quad T_i\mapsto f_iT_i$
are injective (see \cite[4.2 Def. 1]{Bos14}) and since $\oo_{S^{rig}}$ is a sheaf, we have 
\begin{equation*}
\nomenclature{\sortkey{Ba}$B$}{Ring of global sections $\oo_{S^{rig}}(S^{rig})$ for quasi-stein spaces $S^{rig}$ with $S=Spf\ A$}
B:=\oo_{S^{rig}}(S^{rig})=\varprojlim_n \oo_{S^{rig}}(U_n)=\varprojlim_n B_n\otimes_{\fzet}\F_q\semr{\zeta}=\bigcap_{n\in \N}B_n\otimes_{\fzet}\F_q\semr{\zeta}.
\end{equation*} 
The injective maps $A\hookrightarrow B_n\otimes_{\fzet}\fzett$ yield an injective map $A\hookrightarrow B$. Using this morphism, we can view $A\semm$ as a subring of $B\sem{z-\zeta}$ by the morphism
\begin{equation}\label{neub}
\sum_{i\in \infty}^{\infty}a_iz^i\longmapsto \sum_{i\in \infty}^{\infty}a_i(z-\zeta+\zeta)^i
=\sum_{j=0}^\infty\zeta^{-j}\left(\sum_{i=-\infty}^{\infty}
\begin{pmatrix}
i \\ j
\end{pmatrix} a_i\zeta^i \right)(z-\zeta)^j\quad .
\end{equation}
One has to verify that each coefficient $\zeta^{-j}\left(\sum_{i=-\infty}^{\infty}
\begin{pmatrix}
i \\ j
\end{pmatrix} a_i\zeta^i \right)$ lies in $B_n\otimes_{\fzet}\fzett$ for every $n$, but this is ensured by the condition that $\lim_{i\to -\infty}v(a_i)+ki=\infty$ for all $k>0$ and that $B_n$ is $\zeta$-adically complete and separated.

\rem{}{\mbox{}\label{lperiod} We give some remarks on the element $\lper$. First we have $\lper=\left(1-\frac{\zeta}{z}\right)\sigma(\lper)$. 
If we write $\sigma(\lper)=\sum_{i=0}^\infty l_i(z-\zeta)^i$ with $l_i\in\F_q\semr{z}$ then
the constant coefficient $l_0$ of $\sigma(\lper)\in B\sem{z-\zeta}$ is given by $\prod_{i\in\N_{>0}}\left(1-\zeta^{q^i-1}\right)$. 
We obtain it as the "evaluation" of $\sigma(\lper)$ at $z=\zeta$. Since $\prod_{i\in\N_{>0}}\left(1-\zeta^{q^i-1}\right)\in \fzet^{{*}}\subseteq B^{{*}}$, we conclude $\sigma(\lper)\in B\sem{z-\zeta}^{*}$. The above equation $\lper=(z-\zeta)\frac{1}{z}\sigma(\lper)$ shows that $\lper\notin B\sem{z-\zeta}^{*}$, and hence $\lper\notin A\semm^{*}$. But certainly it shows $\lper\in B\sem{z-\zeta}[\frac{1}{z-\zeta}]^{*}$
and therefore we have an injection 
\begin{equation*}
A\semm\left[\frac{1}{\lper}\right]\subseteq B\sem{z-\zeta}\left[\frac{1}{z-\zeta}\right]
\end{equation*}
by the universal property of the localization at the element $\lper$ and the fact that $\lper$ is not a zerodivisor (otherwise $\sigma(\lper)$ would be a zerodivisor).
Furthermore, we note that $A\semm\otimes_A A/I=A/I\semr{z}$ and since $\lper\equiv 1\mod I$ we also have $A\semm\left[\frac{1}{\lper}\right]\otimes_A A/I=A/I\semr{z}$.
}
\nomenclature{\sortkey{lo}$l_0$}{constant coefficient of $\sigma(\lper)$ which is equal to $\prod_{i\in\N_{>0}}\left(1-\zeta^{q^i-1}\right)\in\F_q\semr{z}$}

\lem{{compare also \cite[Lemma 2.3.1]{Har11}
}}{\label{etalim}\mbox{}\\
Let $(M_n, \tau_{M_n}, \eta_n)_{n\in\N}\in\gf(S)$. Then $(\eta_n)_{n\in\N}$ defines a rigid analytic isomorphism 
\begin{equation*}
\eta_M:M_{S^{rig}}\semr{z-\zeta}\to \ds\M\otimes_{\F_q\sem{z},z\mapsto \zeta+(z-\zeta)}\oo_{S^{rig}}\semr{z-\zeta}=:\M_{S^{rig}}\semr{z-\zeta}
\end{equation*}}

Proof: The lemma is proven in \cite[Lemma 2.3.1]{Har11}. We explain the rough idea to make the proof of proposition \ref{etainteger} understandable. After choosing locally a basis we may assume $S=Spf\ A$,  $M\simeq A\semr{z}^r$ as well as $\mathbb M\simeq A\semr{z}^r$. Then $\eta_n$ is given by a matrix in $GL_r(A/I^n\semr{z})$ satisfying $\eta_{n+1}\equiv \eta_n\ mod\ I^n$. The limit $\eta_M:=\varprojlim_n \eta_n$ defines a matrix $Mat_r(\Omega)$, 
where \begin{equation*}
\Omega:=\varprojlim_n \left(A/I^n\semr{z}\right)=\left\{\sum_{i=-\infty}^\infty a_iz^i\ |\ \lim_{i\to \infty}a_i= 0 \right\}.
\end{equation*}
Then one uses in a clever and recursively way the 
equation $\eta_M=\tau_{\M}\circ\sigma^{*}\eta_M \circ \tau_M^{-1}$ to estimate the $I$-valuation of the coefficients of the power series arising in the matrix $\eta_M$. A computation shows then that $\eta_M$ is a matrix in $Gl_r\left(A\semm\left[\frac{1}{l_-}\right]\right)$.\hfill$\Box$

\prop{}{\label{etainteger}\mbox{}
With the previous notation the isomorphism
$$\displaystyle\sigma^{{*}}\eta_M:\sigma^{*} M_{S^{rig}}\semr{z-\zeta}\to \sigma^{*}\M_{S^{rig}}\semr{z-\zeta}$$ restricts to an isomorphism
$$\displaystyle\sigma^{{*}}\eta_M:\sigma^{{*}}M_{S^{rig}}\sem{z-\zeta}\to \sigma^{*}\M_{S^{rig}}\sem{z-\zeta}$$ 
of $\oo_{S^{rig}}\sem{z-\zeta}$-modules.}

\prof{From lemma \ref{etalim} we get an induced isomorphism $\displaystyle\sigma^{{*}}\eta_M:\sigma^{*} M_{S^{rig}}\semr{z-\zeta}\to \sigma^{*}\M_{S^{rig}}\semr{z-\zeta}$.
Let $S=\bigcup_{i\in I}U_i$ be an open affine covering of $S$ with $U_i=Spf\ A_i$ such that $M(U_i)$ is a free $A_i$ module. Then $\bigcup_{i\in I}U_i^{rig}$ gives an admissible open covering of $S^{rig}$ and we set $B_{(i)}=\oo_{S^{rig}}(U_i^{rig})$. In the proof of the previous lemma, we saw that $\eta_M|_{U^{rig}}$ is given by a matrix in $Gl_r\left(A_i\semmm\right)\subseteq Gl_r\left(B_{(i)}\sem{z-\zeta}[\frac{1}{z-\zeta}]\right)$. The fact that $\lper\notin B_{(i)}\sem{z-\zeta}$ from remark \ref{lperiod} on page \pageref{lperiod} shows that the assertion can not generally be true for $\eta_M$ itself. But in remark \ref{lperiod} we also saw that $\sigma (\lper)\in B_{(i)}\sem{z-\zeta}^{*}$ and hence $A_i\semm\left
[\frac{1}{\sigma(\lper)}\right]\subseteq B_{(i)}\sem{z-\zeta}$. 
This implies that $\sigma^\star\eta_M|_{U^{rig}}$ is given by a matrix in $Gl_r\left(A_i\semm\left
[\frac{1}{\sigma(\lper)}\right]\right)\subseteq Gl_r\left(B_{(i)}\sem{z-\zeta}\right)$ and thus it proves the proposition.
}

\mbox{}\\
Using the proposition we can now define the Hodge-Pink lattice associated with $x$. Here we write again $\tau_M$ for the induced morphism $\tau_M^{rig}:\sigma^{\star}{M}\otimes_{\oo_S\lsem z\rsem }\oo_{S^{rig}}\semr{z-\zeta}\to {M}\otimes_{\oo_S\lsem z\rsem }\oo_{S^{rig}}\semr{z-\zeta}$. 

\df{}{\label{defqx} For $x=(M_n, \tau_{M_n}, \eta_n)_{n\in \N}\in \gf(S)$ we define $\q_x$ as the image sheaf of $M_{S^{rig}}\sem{z-\zeta}$ under the morphism $\sigma^{*} \eta_M\circ \tau_M^{-1}:M_{S^{rig}}\semr{z-\zeta}\to \sigma^{*}\M_{S^{rig}}\semr{z-\zeta}$
\begin{equation*}
\q_x:=\sigma^{*} \eta_M\circ \tau_M^{-1}(M_{S^{rig}}\sem{z-\zeta}).
\end{equation*}
}
\noindent It is part of the following lemma, that the definition does not depend on the family $(M_n, \tau_{M_n}, \eta_{n})_{n\in\N}$ representing the point $x\in\gf(S)$.
%
%
%
%
%
%
%
%
\lem{}{\label{gammahodge} For every formal $\F_q\sem{z}$-scheme $S$ locally formally of finite type, we get a well-defined map
\begin{align*}
\gamma(S):&&\gf(S)& \longrightarrow \left\{\minibox{Hodge-Pink lattices $\q$ over $S$ on $\underline\M_{\F_q\semr{\zeta}}$,\\ which are bounded by $\omega=(-\mu_r,\dots,-\mu_1)$ }\right\}\\
&&x& \longmapsto\qquad \q_x
\end{align*}
}
\nomenclature{\sortkey{gamma}$\gamma(S)$}{denotes the map from $\gf(S)$ to Hodge-Pink lattices on a fixed $z$-isocrystal}
\prof{
Let $(M_n,\tau_{M_n},\eta_{M_n})_{n\in\N}$ and $(N_n,\tau_{N_n},\eta_{N_n})_{n\in\N}$ be two families representing the same point $x\in\gf(S)$, which means that there is a family $(g_n)_{n\in \N}$ of isomorphisms $g_n:M_n\to N_n$ with $\eta_{N_n}\circ g_n=\eta_{M_n}$. In the limit, this defines an isomorphism $M\to N$ and $g:M_{S^{rig}}\sem{z-\zeta}\to N_{S^{rig}}\sem{z-\zeta}$ with $\eta_N\circ g=\eta_M$ which is seen from the construction of $\eta_M$, $\eta_N$ and $g$.
Therefore, we get the following commutative diagram:

\centerline{
\xymatrix{
M_{S^{rig}}\semr{z-\zeta}
\ar[d]^{\tau_M^{-1}} \ar[r]^{g}
    &
N_{S^{rig}}\semr{z-\zeta}
\ar[d]^{\tau^{-1}_{N}} \ar[r]^{\eta_{N}}
    &
\M_{S^{rig}}\semr{z-\zeta} \ar[d]^{\tau^{-1}_{\M}\otimes id}
    \\
\sigma^{*} M_{S^{rig}}\semr{z-\zeta}
\ar[r]^{\sigma^{*} g} \ar@/_1.5pc/[rr]_{\sigma^{*}\eta_M}
    &
\sigma^{*} N_{S^{rig}}\semr{z-\zeta}
\ar[r]^{\sigma^{*} \eta_{N}}
    &
\sigma^{*} \M_{S^{rig}}\semr{z-\zeta}
}}
Hence we have the following equation
\begin{equation*}
\q_x=\sigma^{*}\eta_M \circ \tau_M^{-1}(M_{S^{rig}}
\sem{z-\zeta}
)=\sigma^{*}\eta_M \circ \tau_M^{-1}\circ g^{-1}(N_{S^{rig}}\sem{z-\zeta}
)=\sigma^{*}\eta_{N} \circ \tau_{N}^{-1}(N_{S^{rig}}
\sem{z-\zeta}
),
\end{equation*}
which shows that $\q_x$ is well-defined. It remains to show that $\q_x$ is a Hodge-Pink lattice that is bounded by $\omega$.\\
We choose an admissible open affinoid covering $S^{rig}=\bigcup_{i\in I}Sp\ R_i$ such that $M_{S^{rig}}
\sem{z-\zeta}
(Sp\ R_i)$ is a free $R_i\sem{z-\zeta}$ module of rank $r$. After choosing a basis $M_{S^{rig}}
\sem{z-\zeta}
(Sp\ R_i)=R_i\sem{z-\zeta}^r$ the morphism $\tau_M^{-1}|_{Sp\ R_i}$ is given by a matrix in $Gl_r(R_i\semr{z-\zeta})$ also denoted by $\tau_M^{-1}$. In the same way $\sigma^{*} M_{S^{rig}}\sem{z-\zeta}$ and $\sigma^{*} \M_{S^{rig}}\sem{z-\zeta}$ are free $R_i\sem{z-\zeta}$-modules and due to proposition \ref{etainteger}, the morphism $\sigma^{*}\eta_M|_{Sp\ R_i}$ is given by a matrix in $Gl_r(R_i\sem{z-\zeta})$.
We denote the matrix again by $\sigma^{*}\eta_M$. In particular, $\q_x(Sp\ R_i)=\sigma^{*}\eta_M\circ \tau_M^{-1}(R_i\sem{z-\zeta}^r)$ is a free $R_i\sem{z-\zeta}$ submodule of rank $r$ of $\sigma^{*} \M_{S^{rig}}\semr{z-\zeta}(Sp\ R_i)\simeq R_i\sem{z-\zeta}[\frac{1}{z-\zeta}]^r$. This shows that $\q_x$ is a locally free $\oo_{S^{rig}}\sem{z-\zeta}$ module of rank $r$ and that $\q_x\otimes_{\oo_{S^{rig}}\sem{z-\zeta}}\oo_{S^{rig}}\semr{z-\zeta}\simeq \sigma^{*}\M_{S^{rig}}\semr{z-\zeta}$. 
Therefore, $\q_x$ is a Hodge-Pink lattice over $S$ on $\underline\M_{\F_q\semr{\zeta}}$ and we have to prove that it is bounded by $\omega=(-\mu_r,\dots,-\mu_1)$. We recall that $\p_{S^{rig}}=\sigma^{*}\M\otimes_{\F_q\sem{z}}\oo_{S^{rig}}\sem{z-\zeta}$. Since the local Shtuka $\underline M$ is bounded by $\mu=(\mu_1,\dots,\mu_r)$ we have for all  $1\klg i\klg r$:
\begin{alignat*}{4}
&&\quad\bigwedge^i \tau_M\bigwedge^{i}\sigma^{{*}}M_{S^{rig}}\sem{z-\zeta} & \subseteq 
(z-\zeta)^{\mu_{r-i+1}+\dots+\mu_r}\bigwedge^iM_{S^{rig}}\sem{z-\zeta}&\quad &\mbox{}\\
\Leftrightarrow &&\bigwedge^{i}\sigma^{{*}}M_{S^{rig}}\sem{z-\zeta} & \subseteq 
(z-\zeta)^{\mu_{r-i+1}+\dots+\mu_r}\bigwedge^i\tau_M^{-1}(M_{S^{rig}}\sem{z-\zeta})&\quad &
\\
&&&\mbox{(now use $\sigma^{*}\eta_M(\sigma^{*} M_{S^{rig}}\sem{z-\zeta})=\sigma^{*}\M_{S^{rig}}\sem{z-\zeta}$ by \ref{etainteger})}&&\\
\Leftrightarrow &&\bigwedge^{i}\sigma^{{*}}\M_{S^{rig}}\sem{z-\zeta} & \subseteq 
(z-\zeta)^{\mu_{r-i+1}+\dots+\mu_r}\bigwedge^i\sigma^{{*}}\eta_M\circ\tau_M^{-1}(M_{S^{rig}}\sem{z-\zeta})&\quad &
\\
\Leftrightarrow &&\quad(z-\zeta)^{-\mu_r-\dots-\mu_{r-i+1}}\bigwedge^{i}\p_{S^{rig}} & \subseteq 
\bigwedge^i\q_x &\quad &
\end{alignat*}
The equivalences remain true, if we require for the inclusions equality in the case that $i=r$. It follows that $\q_x$ is indeed bounded by $\omega =(-\mu_r,\dots,-\mu_1)$, which proves the lemma.
 }

\bsp{}{\label{beispielhodgelattice}
Let $\gf$ be the Rapoport-Zink functor associated with \casecomp and let 
$\gs=\coprod_{(k,l)\in\Z^2}U_{kl}=\coprod_{\Z^2}Spf\ \F_q\sem{z,h}$ the Rapoport-Zink space and $x^{univ}\in\gf(\gs)$ the universal ind-object as determined in section \ref{universalob}. Since we will need it in the next section, we want to compute the Hodge-Pink lattice $\q_{x^{univ}}=\gamma(\gs)(x^{univ})$. It is sufficient to compute $\q_{x^{univ}}|_{U_{kl}^{rig}}$ restricted to the connected components $U_{kl}^{rig}$. By definition \ref{defqx}, 
it is given as $\q_{x^{univ}}|_{U_{kl}^{rig}}=\sigma^{*}\eta_M\circ \tau_M^{-1}(M_{{RZ_{b,\mu}^{\ rig}}}\sem{z-\zeta}|_{U_{kl}^{rig}})$. Let $B=\oo_{{RZ_{b,\mu}^{\ rig}}}(U_{kl}^{rig})$, then $M_{{\gs}^{rig}}\sem{z-\zeta}|_{U_{kl}^{rig}}$ is isomorphic to the $\oo_{U_{kl}^{rig}}\sem{z-\zeta}$ module associated with $B\sem{z-\zeta}^2$. Using the results from section \ref{universalob} as well as the fact that $\sigma^{*}\eta_M$ and $\tau_M^{-1}$ are given by matrices in $Gl_r(B\sem{z-\zeta}[\frac{1}{z-\zeta}])$, we can write:
\begin{equation*}
M_{{RZ_{b,\mu}^{\ rig}}}\sem{z-\zeta}|_{U_{kl}^{rig}}=B\sem{z-\zeta}^2=B\sem{z-\zeta}\vecc{1}{0}+B\sem{z-\zeta}\vecc{0}{1},
\end{equation*}
\begin{equation*}
\sigma^{*}\eta_M|_{U_{kl}^{rig}}=\varprojlim_n\sigma^{*} \eta_n=
\varprojlim_n
\begin{pmatrix}
z^k\prod_{i=0}^{n-1}\frac{z}{z-\zeta^{q^{i+1}}} & 0 \\
-z^l\sum_{i=0}^{n-1}\frac{{h}^{q^{i+1}}}{(z-\zeta^q)\cdot\dots\cdot(z-\zeta^{q^{i+1}})}& z^l
\end{pmatrix}=
\begin{pmatrix}
z^k\prod_{i=0}^{\infty}\frac{z}{z-\zeta^{q^{i+1}}} & 0 \\
-z^l\sum_{i=0}^{\infty}\frac{{h}^{q^{i+1}}}{(z-\zeta^q)\cdot\dots\cdot(z-\zeta^{q^{i+1}})}& z^l
\end{pmatrix},
\end{equation*}
\begin{equation*} 
\tau_M|_{U_{kl}^{rig}}=
\begin{pmatrix}
z-\zeta & 0 \\ h & 1
\end{pmatrix}
\mbox{ \quad and, therefore, }
\tau_M^{-1}|_{U_{kl}^{rig}}=
\begin{pmatrix}
\frac{1}{z-\zeta} & 0 \\ -\frac{h}{z-\zeta} & 1
\end{pmatrix}.
\end{equation*}
We compute 
\begin{equation*}
\sigma^{*}\eta\circ\tau_M^{-1}|_{U_{kl}^{rig}}=
\begin{pmatrix}
z^{k-1}\prod_{i=0}^{\infty}\frac{z}{z-\zeta^{q^i}} & 0 \\
-z^l\sum_{i=0}^{\infty}\frac{{h}^{q^i}}{(z-\zeta)\cdot\dots\cdot(z-\zeta^{q^i})}& z^l
\end{pmatrix}.
\end{equation*}
With $\frac{1}{z}\prod_{i=0}^{\infty}\frac{z}{z-\zeta^{q^i}}=\frac{1}{z}\lper^{-1}=(z-\zeta)^{-1}\sigma^{*}(\lper^{-1})$ this results in 
\begin{align*}
\q_{x^{univ}}|_{U_{kl}^{rig}}&=\bvecc{\frac{z^k}{z-\zeta}\sigma^{*}(\lper^{-1})}{-z^l\sum_{i=0}^{\infty}\frac{{h}^{q^i}}{(z-\zeta) \dots (z-\zeta^{q^i})}} B\sem{z-\zeta}+\bvecc{0}{z^l}B\sem{z-\zeta}\\
&=\bvecc{1}{-z^{l-k}\sigma^{*}(\lper)\cdot \sum_{i=0}^{\infty}\frac{{h}^{q^i}}{(z-\zeta^q) \dots (z-\zeta^{q^i})}} (z-\zeta)^{-1} B\sem{z-\zeta}+\bvecc{0}{1}B\sem{z-\zeta}.
\end{align*}
}
\section{The Period Morphism}\label{period}

In this section we will define the period morphism $\pi:{RZ_{b,\mu}^{\ rig}}\to \Qsr$, which has the characteristic property that it contains the information about all the maps $\gamma(S):\gf(S)\to \Qf(S^{rig})$ for every formal scheme $S\in \flfts$.

We will also describe how the maps $\gf(S)\to \Qf(S^{rig})$ can be recovered from the period morphism $\pi:{RZ_{b,\mu}^{\ rig}}\to\Qsr$. Furthermore, we will compute explicitly the period morphism in the case where $\underline{\M}=\big(\F_q\sem{z}^2,(\begin{smallmatrix}z&0\\0&1\end{smallmatrix})\big)$ and $\mu=(1,0)$, which turns out to be the Carlitz logarithm.\\

Since $\gf$ and $\Qf$ are representable by $\gs$ and $\Qsr$ the maps $\gamma(S)$ from Lemma \ref{gammahodge} define for every formal scheme $S$ locally formally of finite type a map $\mathfrak{n}(S):Hom(S,\gs)\to Hom(S^{rig},\Qsr)$ by the following diagram:
\begin{equation}\label{funktorendarstellbar}
\xymatrix{
\gf(S)  \ar[d]_{\gamma(S)} \ar[rr]^{\widesim\qquad\qquad} && Hom(S,\gs) \ar[d]^{\mathfrak{n}(S)}\qquad\\
\Qf(S^{rig}) \ar[rr]^{\widesim\qquad\qquad} && Hom(S^{rig},\Qsr)
}
\end{equation}

\df{}{
\label{periodmorphismdef}\mbox{}\\
We define the period morphism $\pi:{RZ_{b,\mu}^{\ rig}}\to \Qsr$ of the moduli problem to be $\mathfrak{n}(\gs)(id_{\gs})\in Hom({RZ_{b,\mu}^{\ rig}},\Qsr)$.
}

The following lemma with its corollary shows that the period morphism has indeed the announced properties.

\lem{}{The maps $\mathfrak{n}(S):Hom(S,\gs)\to Hom(S^{rig},\Qsr)$ define a natural transformation\\\mbox{}\quad $\mathfrak{n}:Hom(-,\gs)\to Hom(-,\Qsr)\circ rig$\quad in the category $Fun\left(\op{(\flfts)},Set\right)$.}

\prof{
The statement of the lemma is equivalent to the assertion that the maps $\gamma(S)$ define a natural transformation $\gamma:\gf\to \Qf\circ rig$. This is clear by the definition of $\gamma(S)$.}

\ko{}{
The natural transformation $\mathfrak{n}:Hom(-,\gs)\to Hom(-,\Qsr)\circ rig$ is uniquely determined by the period morphism $\pi:{RZ_{b,\mu}^{\ rig}}\to \Qsr$.
}
\prof{
This is the statement of the Yoneda lemma (see \cite[3.2 page 61]{Mac71}) applied to the above functors. It tells us that there is a bijection:
\begin{eqnarray*}
nat\big(Hom(-,\gs\big),{Hom(-,\Qsr)\circ rig})&\xlongleftrightarrow{1:1} &{Hom(RZ_{b,\mu}^{\ rig},\Qsr)}\\
\mathfrak{m} & \longmapsto & \mathfrak{m}(\gs)(id_{\gs})\\
\left(\minibox{$\mathfrak{m}(S):Hom(S,\gs)\to{Hom(S^{rig},\Qsr)}$\\
\qquad\quad\qquad$f\mapsto {Hom(-,\Qsr)\circ rig}(f)(X)$}\right) & \ \llmapsto & X\in{Hom(RZ_{b,\mu}^{\ rig},\Qsr)}
\end{eqnarray*}
\noindent This proves the corollary since by definition $\pi=\mathfrak{n}(\gs)(id_{\gs})$.
}\\

The second map in the above proof determines the maps $\mathfrak{n}(S)$ in terms of $\pi$ and the functor ${Hom(-,\Qsr)\circ rig}$. If one knows the universal object $x^{univ}\in\gf(\gs)$ of the representable functor $\gf$ and the Hodge-Pink lattice $\q_{x^{univ}}\in \Qf({RZ_{b,\mu}^{\ rig}})$, one can describe the maps $\gamma(S)$ explicitly in terms of $\gf,\Qf, \pi, x^{univ}$ and $\q_{x^{univ}}$.

\noindent Namely let $x\in \gf(S)$ with $x=\gf(f)(x^{univ})$ for a unique morphism $f:S\to \gs$. Then equation \eqref{funktorendarstellbar} suggests:
\begin{equation*}
\gamma(S)(x)=\Qf(\pi\circ f^{rig})(\q_{x^{univ}})
\end{equation*}
\noindent More concretely, we can also write $\gamma(S)(x)=\q_{x^{univ}}\otimes_{\oo_{{RZ_{b,\mu}^{\ rig}}}\sem{z-\zeta}}\oo_{S^{rig}}\sem{z-\zeta}$, where the morphism ${\oo_{{RZ_{b,\mu}^{\ rig}}}\sem{z-\zeta}}\to \oo_{S^{rig}}\sem{z-\zeta}$ used in the tensor product is induced by the morphism $\pi\circ f^{rig}=\mathfrak{n}(S)(f)\in Hom(S^{rig},\Qsr)$.

 \begin{center}
{ \textbf{The Carlitz Logarithm}}
 \end{center}
 
We recall the definition of the Carlitz logarithm. Let $\delta$ be the morphism $\delta:\F_q[t]\to \F_q(\zeta)\quad t\mapsto \zeta$. Moreover we fix the valuation $|-|:\F_q(\zeta)\to \R^+_0$ with $|\zeta|=\frac{1}{q}$. Then $\F_q\semr{\zeta}$ is the completion with respect to this valuation and we denote by $\C_t=\widehat{\F_q\semr{\zeta}^{alg}}$ the completion of its algebraic closure. 
In addition we write $B(0,\rho)=\{x\in \C_t\ \big|\ |x|<\rho\}$.
As usual, let $\F_q(\zeta)\{\tau\}$ (resp. $\F_q(\zeta)\{\{\tau\}\}$) be the noncommutative polynomial ring (power series) with $\tau b=b^q\tau$ for $b\in\F_q(\zeta)$. The Carlitz module is defined as $\varphi:\F_q[t]\to \F_q(\zeta)\{\tau\}$ with $t\mapsto \zeta+\tau$. We will write $\varphi_a$ instead of $\varphi(a)$. Viewing $\tau$ as a morphism $\tau:\F_q(\zeta)\to \F_q(\zeta)$ with $b\mapsto b^q$, it defines a $\F_q[t]$-module structure on $\F_q(\zeta)$ by $a\cdot b=\varphi_a(b)$ for $a\in\F_q[t]$ and $b\in\F_q(\zeta)$.

\df{{\cite[Definition 3.2.7]{Gos91}}}{\mbox{}
The Carlitz exponential is defined to be the formal power series
\begin{equation*}
\exp_\varphi=\sum_{n=0}^\infty \frac{\tau^n}{(\zeta^{q^n}-\zeta)\dots(\zeta^{q^n}-\zeta^{q^{n-1}})}\in \F_q(\zeta)\{\{\tau\}\}
\end{equation*}
}

This defines a map $\exp_\varphi:B(0,|\zeta|^{1-q})\to \C_t,\quad x\mapsto \sum_{n=0}^\infty e_n x^{q^n}$
with 
$e_n=\frac{1}{(\zeta^{q^n}-\zeta)\dots(\zeta^{q^n}-\zeta^{q^{n-1}})}$.
Since we have $|e_n|=|\zeta|^{-q^{n-1}...-q-1}$ and therefore $\sqrt[q^n]{|e_n|}=|\zeta|^{-q^{-1}-q^{-2}\dots -q^{-n}}$ the radius of convergence is:
\begin{equation*}
r(\exp_\varphi)=(\limsup_{n\to \infty}\sqrt[q^n]{|e_n|})^{-1}=(|\zeta|^{-\frac{1}{q}\frac{1}{1-q^{-1}}})^{-1}=|\zeta|^{1-q}
\end{equation*}

We remark that the Carlitz exponential satisfies\begin{equation*}
\exp_\varphi \cdot D_0(\varphi_a) = \varphi_a(\exp_\varphi)\mbox{ for all $a\in\F_q[t]$, where  $D_0(\varphi_a):=a_0$ for $\varphi_a=\sum_{i=0}^m a_i\tau^i$}.
\end{equation*}

\df{}{We define the Carlitz logarithm as 
\begin{equation*}
\log_\varphi:=\exp_\varphi^{-1}\in \F_q(\zeta)\{\{\tau\}\}
\end{equation*}
}

The previous remark yields $D_0\varphi_a\cdot \log_\varphi=\log_\varphi\cdot \varphi_a$ and in particular $\zeta \cdot \log_\varphi=\log_\varphi \cdot \varphi_t$.
If we write $\log_\varphi=\sum_{i=0}^\infty c_n\tau^n$ this implies

\begin{equation*}
\sum_{n=0}^\infty\zeta c_n\tau^n = \sum_{n=0}^\infty(c_n\zeta^{q^n}+c_{n-1})\tau^n\quad \Leftrightarrow\ (\zeta^{q^n}-\zeta)c_n=-c_{n-1},\quad\mbox{ for } n>0\mbox{ and }c_0=1
\end{equation*}

It follows that $c_n=\prod_{i=1}^n\frac{1}{(\zeta-\zeta^{q^i})}\in \F_q(\zeta)$
and for the absolute value of the coefficients we get 
$|c_n|=|\zeta|^{-n}$ and, therefore, $r(\log_\varphi)=(\limsup_{n\to\infty}\sqrt[q^n]{|c_n|})^{-1}=1$. We conclude:

\lem{}{\label{carlitzlogdef}
The Carlitz logarithm $\log_\varphi$ is given by 
\begin{equation*}
\log_\varphi(X):=\sum_{n=0}^\infty\frac{X^{q^n}}{(\zeta-\zeta^{q})\dots(\zeta-\zeta^{q^n})}\in \F_q(\zeta)\lsem X\rsem\subset\C_{(t)}\lsem X\rsem
\end{equation*}
and has the radius of convergence $1$.}

\begin{center}\textbf{
	The Period Morphism for the Moduli Problem 
  in section \ref{particular}}
\end{center}

In the final part of this article we explicitly compute the period morphism for the Rapoport-Zink functor
$\gf$ associated with 
 $\underline\M=\big(\F_q\sem{z}^2,(\begin{smallmatrix}
z&0\\0&1
\end{smallmatrix}
) \big)$ and $\mu=(1,0)$.
However, throughout this article we have been building towards this explicit computation.
We proved that $\gf$ is representable by the formal $\F_q\sem{z}$-scheme $\gs=\ds\coprod_{(i,j)\in\Z^2}Spf\ \F_q\sem{\zeta,h}$, we described 
the universal ind-object $x^{univ}\in \gf(\gs)$, we computed the Hodge-Pink lattice $\q_{x^{univ}}=\gamma(RZ_{b,\mu})(x^{univ})$ and we showed that the functor 
$\Qf=\Qf_{D,\omega}$ with $\underline D=\big(\F_q\semr{z}^2,(\begin{smallmatrix}
z&0\\0&1
\end{smallmatrix}
) \big)$ and $\omega=(0,-1)$ is represented by $\PP_{\fzett}^{1,rig}=Grass^{rig}_{1,2}$.

The period morphism $\pi$ was defined as $\mathfrak{n}(\gs)(id_{\gs})$ and it is also given as the image of  $\q_{x^{univ}}$ under the lower map in diagram \eqref{funktorendarstellbar}. It suffices again to compute $\pi$ on each connected component $U_{kl}$.

\noindent In example \ref{beispielhodgelattice} we computed $\q_{x^{univ}}|_\uklr$ as the $\oo_{{U^{rig}_{k,l}}}$ module associated with the module
\begin{equation*}
\bvecc{(z-\zeta)^{-1}}{\ds-(z-\zeta)^{-1}\sigma^{*} \lper\cdot z^{l-k} \sum_{i=0}^\infty\frac{h^{q^i}}{(z-\zeta^q)\dots(z-\zeta^{q^i})}}B\sem{z-\zeta}+\bvecc{0}{1}B\sem{z-\zeta}\subseteq B\sem{z-\zeta}\left[\frac{1}{z-\zeta}\right]^2
\end{equation*}
where $B=\oo_{{RZ_{b,\mu}^{\ rig}}}(\uklr)$. In example \ref{beispielq}, we saw that in

 our particular case $\Qf$ is isomorphic to $Grass^{rig}_{1,2}$ by the map 
\begin{equation*}
\Qf(S)\isom Grass^{rig}_{1,2}\qquad \q_x\mapsto Fil^1\ D_{S^{rig}}(\q_x).
\end{equation*}

\noindent Recall from remark \ref{lperiod} on page \pageref{lperiod} that $l_0:=\prod_{i\in\N_{>0}}\left(1-\zeta^{q^i-1}\right)\in \fzett$ is the constant coefficient of $\sigma(\lper)\in B\sem{z-\zeta}$.
 We have $D_{\uklr}=\sigma^{*} \M_\uklr[\frac{1}{z}]\simeq \oo_\uklr^2$ and $Fil^1\oo_\uklr^2(\q_{x^{univ}})$ is the locally free $\oo_\uklr$-submodule of $\oo_\uklr^2$ associated with the module 
\begin{alignat*}{2}
&\bruch{B\sem{z-\zeta}^2\cap(z-\zeta)\q_{x^{univ}}|_{\uklr}}{(z-\zeta)B\sem{z-\zeta}^2\cap (z-\zeta)\q_{x^{univ}}|_\uklr}\\
%
%
=\quad & \bruch{
\bvecc{1}{
\ds-\sigma^{*} \lper\cdot z^{l-k} \sum_{i=0}^\infty\frac{h^{q^i}}{(z-\zeta^q)\dots(z-\zeta^{q^i})}
}B\sem{z-\zeta}
{
+\bvecc{0}{
z-\zeta
}B\sem{z-\zeta}
}
}
{(z-\zeta)B\sem{z-\zeta}^2}\\
%
%
%
 =\quad & \bvecc{1}{
-l_0\ds\cdot \zeta^{l-k} \sum_{i=0}^\infty\frac{h^{q^i}}{(\zeta-\zeta^q)\dots(\zeta-\zeta^{q^i})}
}B\quad .
\end{alignat*}

\noindent By lemma \ref{carlitzlogdef} we have $\ds\log_\varphi(h)=\sum_{i=0}^\infty\frac{h^{q^i}}{(\zeta-\zeta^q)\dots(\zeta-\zeta^{q^i})}$.
Next, we just have to determine the image of $Fil^1\oo_\uklr^2(\q_{x^{univ}})$ under the map
$Grass_{1,2}^{rig}(S^{rig})\to \Hom(S^{rig},\PP_{\fzett}^{1,rig})$. 
Since 
\begin{equation*}
\bruch{B^2}{\bvecc{1}{-l_0\ds\cdot \zeta^{l-k} 
\log_\varphi(h)
}B}
\simeq B,\qquad \overline{\bvecc{a}{b}}\mapsto b+a\cdot l_0\ds\cdot \zeta^{l-k} 
\log_\varphi(h)
\end{equation*}
we have the following exact sequence
\begin{alignat*}{3}
0\longrightarrow  \bvecc{1}{-l_0\ds\cdot \zeta^{l-k} 
\log_\varphi(h)
}B \longrightarrow\quad  & B^2\ & \xrightarrow{s_1\oplus s_2} &\quad B &\longrightarrow  0\\
&\bvecc{a}{b} & \longmapsto\ & b+a\cdot l_0\ds\cdot \zeta^{l-k} 
\log_\varphi(h)&\qquad \quad .
\end{alignat*}
This corresponds to an exact sequence of $\oo_\uklr$-modules

\begin{alignat*}{3}
0\longrightarrow  Fil^1 \oo_\uklr^2(\q_x^{univ})
 \longrightarrow\quad  & \oo_\uklr^2\ & \xrightarrow{s_1\oplus s_2} &\quad \oo_\uklr &\longrightarrow  0\\
&\bvecc{a}{b} & \longmapsto\quad & s_1(a)+s_2(b)
\end{alignat*}
where $s_1$ and $s_2$ are module homomorphisms in $\Hom(\oo_\uklr,\oo_\uklr)$, with $s_1=id$ and $s_2$ the homomorphism associated with the multiplication by $l_0\ds\cdot \zeta^{l-k} \log_\varphi(h)$ in $B$. 
We cover $\PP_{\fzett}^{1}=V_1\cup V_2$ with $V_1=Spec\ \fzett [\frac{z_0}{z_1}]$ and $V_2=Spec\ \fzett [\frac{z_1}{z_0}]$.
Then $\PP_{\fzett}^{1,rig}=V_1^{rig}\cup V_2^{rig}$ is covered by $V_1^{rig}$ and $V_2^{rig}$.
Now the map 
\begin{equation*}
f_{k,l}:\fzett\left[\frac{z_0}{z_1}\right]\to B,\qquad \frac{z_0}{z_1}\mapsto \frac{s_2|_{\uklr}(1)}{s_1|_{\uklr}(1)}=\frac{l_0\cdot \zeta^{l-k}\log_\varphi(h)}{1}
\end{equation*}
induces a morphism $Sp\big(B_n\otimes_{\fzet}\fzett\big)\to V_1$ for every $n\in\N$
and therefore a morphism $\uklr\to V_1$. And by the universal property of the rigid analytification 
this defines a unique morphism 
$\uklr\to V_1^{rig}\subseteq\PP_{\fzett}^{1,rig}$.
As it is the image of $\q_{x^{univ}}|_{\uklr}$ in the diagram \eqref{funktorendarstellbar}, this morphism is now the desired period morphism $\pi$ restricted to $\uklr$. 
The global period morphism $\pi$ is then obtained by gluing the morphisms $\pi|_\uklr$ together. 

In particular, we see that $\pi$ factors through the admissible open $V_1^{rig}$ and that it is given up to the constant $l_0\cdot \zeta^{l-k}$ by the Carlitz logarithm.

\addcontentsline{toc}{section}{References}

\printbibliography

\vfill

Paul Breutmann, CNRS, paul.breutmann@imj-prg.fr\\
Institut de Math\'ematiques de Jussieu - Paris Rive Gauche (IMJ-PRG)\\
4 place Jussieu,\ 75252 Paris Cedex 5,\ France

\end{document}